\documentclass{amsart}
\usepackage{hyperref}

\usepackage[english]{babel} \usepackage{graphicx}

\usepackage{amssymb}
\usepackage{amsfonts}
\usepackage{amscd}

\newtheorem{theorem}{Theorem}[section]

\theoremstyle{definition}
\theoremstyle{remark}
\numberwithin{equation}{section}

\newcommand{\Tang}{\mathrm{Tang}}
\newcommand{\Tan}{\mathrm{Tan}}
\newcommand{\Cont}{\mathrm{Cont}}
\newcommand{\tang}{\mathrm{tang}}
\newcommand{\Ls}{\mathrm{Ls}}
\newcommand{\Li}{\mathrm{Li}}
\newcommand{\dist}{\mathrm{d}}
\newcommand{\cl}{\mathrm{cl}}
\newcommand{\intt}{\mathrm{int}}

\begin{document}
\title[Necessary optimality conditions]{Towards historical roots of
necessary conditions of optimality: Regula of Peano}
\author{Szymon Dolecki}
\address{Mathematical Institute of Burgundy\\
Burgundy University, B.P. 47 870, 21078 Dijon, France}
\email{dolecki@u-bourgogne.fr}
\thanks{The first author gratefully acknowledges a partial support by
Dipartimento di Matematica, Universit\`{a} di Trento.}
\author{Gabriele H. Greco}
\address{Dipartimento di Matematica\\
Universit\`{a} di Trento, 38050 Povo (TN), Italy}
\email{greco@science.unitn.it}
\dedicatory{Dedicated to professor Stefan Rolewicz on his 75-th birthday\\and\\Commemorating the 150-th anniversary of the birth of Giuseppe Peano}
\date{%
July 17, 2007
}

\begin{abstract}
At the end of 19th century Peano discerned vector spaces, differentiability,
convex sets, limits of families of sets, tangent cones, and many other
concepts, in a modern perfect form. He applied these notions to solve
numerous problems. The theorem on necessary conditions of optimality (\emph{%
Regula}) is one of these. The formal language of logic that he developed,
enabled him to perceive mathematics with great precision and depth. Actually
he built mathematics axiomatically based exclusively on logical and
set-theoretic primitive terms and properties, which was a revolutionary
turning point in the development of mathematics.
\end{abstract}

\maketitle

\section{Introduction}

The aim of this paper is to trace back the evolution of mathematical
concepts in the work of Giuseppe Peano (1858-1932) that are constituents of 
\emph{Regula}, that is, Peano's theorem on necessary conditions of optimality.

Well-known necessary conditions of maximality of a function at a point, are
formulated in terms of derivative of the function and of tangent cone of the
constraint at that point. Consider a real-valued function $f:X\rightarrow 
\mathbb{R}$, where $X$ is a Euclidean affine space, and a subset $A$ of $X$%
.\bigskip

\textbf{Regula} (of Optimality)\label{regula} \emph{If }$f$\emph{\ is
differentiable at }$x\in A$\emph{\ and }$f(x)=\max \{f(y):y\in A\}$\emph{,
then}%
\begin{equation}
\langle Df(x),y-x\rangle \leq 0\emph{\ for\ every\ }y\in \Tang(A,x)%
\emph{.}  \label{Pp}
\end{equation}%
\smallskip

The \emph{derivative}  $Df(x)$  is defined to be \emph{the} vector $Df(x)$ such that%
\begin{equation}
\lim\nolimits_{y\rightarrow x}\dfrac{f(y)-f(x)-\langle Df(x),y-x\rangle }{%
|y-x|}=0.\label{frechet_diff}
\end{equation}

The \emph{affine tangent cone} $\Tang(A,x)$ of $A$ at $x$ (for
arbitrary $x\in X$) is given by%
\begin{equation}
\Tang(A,x):=\Ls_{\lambda \rightarrow +\infty
}\left( x+\lambda (A-x)\right) ,  \label{tang+def}
\end{equation}%
where the\emph{\ upper limit} $\Ls_{\lambda \rightarrow
+\infty }A_{\lambda }$ of sets $A_{\lambda }$ (as $\lambda $ tends to $+\infty 
$) is defined by%
\begin{equation}
\Ls_{\lambda \rightarrow +\infty }A_{\lambda }:=\{y\in
X:\liminf_{\lambda \rightarrow +\infty }\dist(y,A_{\lambda
})=0\}.  \label{upper+limit}
\end{equation}%
Of course, $\rm{d}$ is the distance and in (\ref{tang+def}), $A_{\lambda }:=x+\lambda (A-x):=\{x+\lambda
(a-x):a\in A\}$.

It is generally admitted among those who study optimization, that modern
definition of differentiability was introduced by Fr\'{e}chet \cite[%
(1911)]{Frechet_diff}, of tangent cone by Bouligand \cite[(1932)]%
{bouligand}, and of limit of sets by Painlev\'{e} \cite[(1905) p.\,8]{zoretti_J} \footnote{%
  A Painlev\'{e}'s student, Zoretti
(1880-1948), attests in \cite[(1912) p.\,145]{zoretti} that Painlev\'{e}
introduced   both upper and lower limit of a family of sets.  Following Zoretti,  Hausdorff in \cite[(1927) p.280]{hausdorff_mengenlehre} and  Kuratowski  in \cite[(1928) p.\,169]{kuratowski28} reiterate this attribution to Painlev\'e.  More clearly, Zoretti calls ``\emph{set-limit}'' of a sequence of sets  the today's upper limit; while by ``\emph{point-limit}'' he means   ``point belonging to the lower limit'' of the sequence.  Painlev\'e's use of  the notion of  ``\emph{set-limit}''
is dated 1902 in \cite[(1905) p.\,8]{zoretti_J}; on the other hand, one finds in \cite[(1909) p.\,8]{zoretti_B} the first occurrence of the notion of   ``\emph{point-limit}'' without any reference to Painlev\'e. Both ``\emph{set-limit}'' and ``\emph{point-limit}'' are present also in \cite[(1912) p.\,145]{zoretti}. }. So we were very surprised to discover that \emph{%
Regula} was already known by Peano in 1887. Indeed, in order to formulate
it, one needed to possess the notions of differentiability and of affine
tangent cone, hence also that of limit of sets.

But our surprise was even greater, because not only all these notions were
familiar to Peano at the end of 19-th century, but they were formulated in a
rigorous, mature way of today mathematics, in contrast with the approximated
imprecise style that dominated in mathematical writings in those times, and
often persisted during several next decades.

Impressed by the so early emergence of these notions, we started to peruse
the work of Peano in order to understand the evolution of the ideas that
lead to \emph{Regula}.

All the citations of Peano's activity, related to the concepts involved in
the optimality conditions, prevalently concern the span of time between the
first appearance of \emph{Regula} in \textquotedblleft \textbf{Applicazioni
Geometriche}\textquotedblright\ \cite[(1887)]{peano87} and its ultimate form
in \textquotedblleft \textbf{Formulario Mathematico}\textquotedblright%
\footnote{%
The previous four editions of \emph{Formulario mathematico} are \emph{%
Formulaire Math\'ematique} tome 1 (1895), tome 2 (1899), tome 3 (1901) and
tome 4 (1903). The first half of the fifth edition was printed in 1905; the
other half earlier in 1908. The ``Index and Vocabulary'' to \emph{Formulario
mathematico} of 1908 was published separately in 1906.} \cite[(1908) p.\,335]%
{peano_1908}, where \emph{Regula} is stated exactly as above.

Tracing back the development and applications of differentiability,
tangency, limit and other concepts, in the work of Peano over the years, we
see evolution and enrichment of their facets. Peano built mathematics
axiomatically, based exclusively on logical and set-theoretic primitive
terms and properties. This was a revolutionary turning point in the
development of mathematics. The reduction of every mathematical object to
the founding concept of \textquotedblleft set\textquotedblright\ (\emph{%
genus supremum}) of Cantor, enabled the emergence of new concepts
related to properties of sets, unconceivable otherwise. Early and
illuminating examples of the fecundity of Cantor's views are in the books 
\emph{Fondamenti per la teorica delle funzioni di variabili reali} of 
Dini \cite[(1878)]{dini78}, \emph{Calcolo differenziale e integrale} of
 Genocchi and Peano \cite[(1884)]{Genocchi}, and the second edition of 
\emph{Cours d'Analyse} of  Jordan \cite[(1893-96)]{Jordan2}.

To appreciate the novelty of Cantor's approach to mathematics, we should
remember the opposition of some luminaries of mathematics that existed at
the beginning of the twentieth century. For example, in the address to the
Congress of Mathematicians in Rome \cite[p.\,182]{poincare} in 1908, 
Poincar\'{e} said

\begin{quotation}
Quel que soit le rem\`{e}de adopt\'{e} [contre le \textquotedblleft
cantorisme\textquotedblright ], nous pouvons nous promettre la joie du m\'{e}%
decin appel\'{e} \`{a} suivre un beau cas pathologique.
\end{quotation}

This paper is not a definitive word on historical roots of conditions of
optimality. For instance, a confront of (\ref{Pp}) with virtual work
principle has still to be addressed\footnote{%
If a force acting on a material point in equilibrium $x$, has a potential $f$%
, then the virtual work principle states 
\begin{equation*}
\delta L:=\langle Df(x),\delta x\rangle \leq 0\text{ for each }\delta x,
\end{equation*}%
where $\delta x$ is a virtual displacement of that point with respect to an 
\emph{ideal constraint} $A$ (either bilateral or unilateral) independent of
the time (see, for example, Banach \cite[(1951)]{Banach_mech}).}. We
have found no evidence of this relationship in the work of Peano, but it is
plausible that he was aware of it (remark that \emph{Regula} is placed in 
\emph{Formulario Mathematico} within the context of mechanics).

This article concerns several historical aspects. From a methodological
point of view, we are focused on primary sources, and not on secondary
founts, that is, on mathematical facts, and not on opinions or
interpretations of other scholars of history of mathematics. On the other
hand, we will avoid to mention, if not necessary, historical facts that are
well-known among those who study optimization (see, for example, Rockafellar
and Wets \cite{rock-wets}, Borwein and Lewis \cite{borwein-lewis}, Aubin 
\cite{aubin}, Aubin and Frankowska \cite{aubin-frank}, Hirriart-Urruty and
Lemar\'{e}chal \cite{Hirriart}, Pallaschke and Rolewicz \cite{rolewicz}).

\section{Affine and vector spaces}

\emph{Applicazioni Geometriche} is based on the extension theory (\emph{Ausdenungslehre}, 1844 edition) of 
Grassmann, presented in detail in \textquotedblleft
Calcolo Geometrico\textquotedblright\ \cite[(1888)]{calcolo88}, where,
forgoing the philosophical aura founding the work of Grassmann, Peano
introduces the modern notion of vector space.

In Grassmann's work, points and vectors coexist distinctly in a common
structure, together with other objects, like exterior products of points and
vectors (see Greco and Pagani \cite{greco_grass}). This subtle distinction was very
demanding in comparison with today habits of mathematicians. Peano maintains
the distinction. For instance, a difference $y-x$ is a vector if both $y$
and $x$ are either points or vectors; otherwise, it is a point (if $x$ is a
vector and $y$ is a point) or a point of mass $-1$ (if $y$ is vector and $x$ is a point).

Moreover, Peano follows Grassmann in construction of metric concepts from
the scalar product of vectors (introduced by Grassmann in \cite[(1847)]{grass_G}).
Following Grassmann and Hamilton, he conceives the gradient of a function as
a vector, differently from a common habit (of using the norm of the
gradient) that prevailed at the pre-vectorial epoch \footnote{%
These observations are relevant for the understanding of Peano's
interpretation of the formula $\langle Df(x),y-x\rangle \leq 0$ that appears
in \emph{Regula}.}.

In several papers Peano applies the geometric calculus of Grassmann, for
instance, to define area of a surface (see \cite[(1887) p.\,164]{peano87}
and \cite[(1890)]{peano_area}) and to give in \cite[(1898)]{peano_vettori} an axiomatic refoundation (today
standard) of Euclidean geometry, based on the primitive notions of point,
vector and scalar product.

Peano's approach to the definition of linear map was slightly different from
(but equivalent to) that commonly adopted nowadays. Peano says that a map $g$
between spaces is \emph{linear} if is additive and bounded, that is, $%
g(x+y)=g(x)+g(y)$ for all $x$ and $y$, and if $\sup \{\left\vert
g(x)\right\vert :|x|<1\}$ is finite. The reader has certainly observed that
the today condition of \emph{homogeneity} is substituted by that of \emph{%
boundedness }\footnote{%
In other moments (for example, in fourth edition of \emph{Formulaire Math%
\'{e}matique} \cite[(1903) p.\,203]{peano_F4}) Peano adopts a different (but
equivalent) definition of linearity, replacing \emph{boundedness} with \emph{%
continuity}. All these variants are related to the following fundamental
lemma (see a proof in \emph{Formulario Mathematico} \cite[(1908) pp.\,117-118]%
{peano_1908}, where Peano quotes Darboux \cite[footnote of p.\,56]%
{darboux_lem}): \textquotedblleft \emph{For an additive function $f:\mathbb{R%
}^{n}\rightarrow \mathbb{R}$ are equivalent: (1) homogeneity, (2) continuity
and (3) boundedness on bounded sets}\textquotedblright}. For Peano, the
interest of employing boundedness in the definition, was to obtain
simultaneously a concept of \emph{norm} (\emph{module} in his terminology)
on spaces of linear maps.

The norm was useful in his study of systems of linear differential equations 
\cite[(1888)]{peano_ED}; to give a formula for a solution in terms of
resolvent, he defines the exponential of matrix and proves its convergence
using the norm (see also Peano \cite[(1894)]{peano_exp} and the English
translation of \cite[(1888)]{peano_ED} in \cite{birkhoff} by G.
Birkhoff).

As other new theories, the theory of vector spaces was contested by many
prominent mathematicians. Even those (few) who adopted the vector approach,
were not always entirely acquainted with its achievements. To perceive the
atmosphere of that time, we give an excerpt from the introduction of 
Goursat to book \emph{Le\c{c}ons de g\'{e}om\'{e}trie vectorielle} \cite[%
(1924)]{Boul_Vect} of Bouligand:

\begin{quotation}
Si le calcul vectoriel a \'{e}t\'{e} un peu lent \`{a} p\'{e}n\'{e}trer en
France, il est bien certain que la multiplicit\'{e} des notations et l'abus
du symbo\-lisme ont justifi\'{e} en partie la d\'{e}fiance de nos \'{e}%
tudiants. Or, dans le livre de M. Bouligand, le symbolisme est r\'{e}duit au
minimum, et l'auteur n'h\'{e}site pas \`{a} revenir aux proc\'{e}d\'{e}s
habituels du calcul quand les m\'{e}thodes lui paraissent plus directes.
[...]

M. Bouligand a devis\'e son ouvrage en trois parties, consacr\'ees
respectivement aux op\'erations vectorielles en g\'eom\'etrie lin\'eaire, en
g\'eom\'etrie m\'etrique et aux op\'erations infinit\'esimales.
\end{quotation}

A decisive role in the dissemination of vector spaces had a book {\em Space-Time-Matter} \cite[(1918)%
]{Weyl} of  Weyl (for details see Zaddach \cite[(1988)]{zaddach88}, \cite[(1994]{zaddach94}).

\section{Differentiability}

In \emph{Applicazioni Geometriche} (p.\,131) Peano says that a
vector $\mathbf{u}$ is a \emph{derivative} at a point $x$ of a real-valued
function $f$ defined on a finite-dimensional Euclidean affine space $X$, if there
exists a vector $\varepsilon (y)$ such that 
\begin{equation}
f(y)-f(x)=\langle y-x,\mathbf{u}+\varepsilon (y)\rangle \text{ with }%
\lim\nolimits_{y\rightarrow x}\varepsilon (y)=0.  \label{differ}
\end{equation}%
The reader recognizes in (\ref{differ}) the Taylor formula of order $1$ and,
on the other hand, the characterization of derivability which 
Carath\'{e}odory gives in \cite[p.\,119]{caratheodory}: \textquotedblleft $f$
is derivable at $x$ if there is a function $\varphi $ continuous at $x$ such
that $f(y)-f(x)=\langle y-x,\varphi (y)\rangle $ for every $y$%
\textquotedblright\ \footnote{%
Remarkably, this Carath\'{e}odory reformulation \textquotedblleft leads to
some sharp, concise proofs of important theorems: chain rule, inverse
function theorem, ...\textquotedblright\ (see Kuhn \cite[(1991)]{kuhn}) and
\textquotedblleft makes perfect sense in general linear topological
spaces\textquotedblright\ (see Acosta-Delgado \cite[(1994)]{acosta}).}.

In \emph{Formulario Mathematico} of 1908 (p.\,334 and 330) the derivative $\mathbf{u}
$ is denoted by $Df(x)$ and is defined by (\ref{frechet_diff}) and, more generally,  one finds  a definition of
differential of map between finite-dimensional Euclidean vector spaces, namely if $f:$ $\mathbb{R}%
^{m}\rightarrow \mathbb{R}^{n}$ then a \emph{derivative} of $f$ at $x$ is 
\emph{the} linear map $Df(x):\mathbb{R}^{m}\rightarrow \mathbb{R}^{n}$
(referred to as \emph{Jacobi-Grassmann} derivative by Peano \cite[(1908) p.\,455]{peano_1908}, and called  nowadays the \emph{Fr\'{e}chet derivative} of $f$ at $x$) such that%
\begin{equation}
\lim\nolimits_{y\rightarrow x}\dfrac{f(y)-f(x)-Df(x)(y-x)}{|y-x|}=0.
\label{der+vect}
\end{equation}

In giving this definition, Peano refers to the second edition of \emph{%
Ausdehnungslehre} of 1962 of Grassmann \cite[v.\,2 p.\,295]{grass} and to an
article \emph{De determinantibus functionalibus} 1841 of  Jacobi \cite[%
v.\,3, p.\,421]{jacobi}. Actually the citation of Jacobi refers to the
concept of \emph{Jacobian}.

With respect to Peano's quotation of Grassmann, our verification of the
source leads the following facts. For a map $f:\mathbb{R}^{m}\rightarrow 
\mathbb{R}^{n}$, Grassmann defines a \emph{differential} $df(x)$ at $x$, as
a map from $\mathbb{R}^{m}$ to $\mathbb{R}^{n}$, by 
\begin{equation}
df(x)(v):=\lim\nolimits_{q\rightarrow 0}\dfrac{f(x+q\,v)-f(x)}{q}.
\label{der+grass}
\end{equation}%
Grassmann proves that if the differential $df(x)$ exists at every $x$ and it is  radially continuous \footnote{%
A map $g$ is \emph{radially continuous at }$x$ if for every vector $v$ the
map $h\longmapsto g(x+hv)$ is continuous at $0$. Grassmann adopts the term 
\emph{continuous} to denote ``radially continuous''.} in variable $x$ for each fixed $v\in V$, then  $f$ is radially
continuous, the differential is linear in $v$, and that (\ref{der+grass})
becomes the partial derivative when $v$ is an element of the canonical base.
Moreover, he claims that the chain rule holds. In contrast to the comments 
\cite[(2000) p.\,398]{grass_E} of Kannenberg (author of English translation \cite%
{grass_E} of \emph{Ausdehnungslehre} 1962), this claim is however false, as
can be seen from Acker-Dickstein's example \cite[Ex.\,2.5, p.\,26]{acker} 
\footnote{%
A counterexample to Grassmann's claim is given by the functions $f$, and $g$
of two real variables defined by $g(x,y):=(x,y^{2})$, $f(x,y):=\frac{x^{3}y}{%
x^{4}+y^{2}}$ and $f(0,0):=0$. These functions was used in \cite{acker} to
invalidate a similar claim for Gateaux differentiability.}.

In various moments of his activity, Peano studied the concept of derivative.
His contributions are very rich and diversified, and concern

\begin{enumerate}
\item strict derivative (see Peano \cite[(1884)]{peano_jordan}, \cite[(1884)]{peano_gilbert} and 
\cite[(1892)]{peano_mathesis})\footnote{%
In \cite[(1884)]{peano_jordan} Peano observes the equivalence between continuity of
derivative  $f'$  of $f$ at $x$ and ``strict derivability of $f$ at $x$'' (that is, $\lim_{a,b\to x}\frac{f(b)-f(a)}{b-a}=f'(x)$). In \cite[(1884)]{peano_gilbert} he notices that
the uniform convergence of the difference quotient function $\frac{f(x+h)-f(x)}{h}%
$ in variable $x$ to $f^{\prime }(x)$  (as h tends to 0) amounts to continuity of
derivative $f^{\prime }$ in variable $x$. As observed Mawhin in \cite[(1997) p.\,430]{mawhin}, Peano
formulates in \cite[(1884)]{peano_gilbert} an approximation property of primitives
equivalent to Kurzweil integrability of all the functions having a primitive.%
},

\item Taylor formula with infinitesimal remainder (called \emph{Peano
remainder}) (see Genocchi \cite[(1884) p.\,XIX]{Genocchi} or Stolz \cite[(1893), pp.\,90-91]{stolz}%
),

\item asymptotic development and Taylor formula in Peano \cite[(1891)]{peano_taylor} 
\footnote{%
If $f$ is a function, and $P(h)=a_{0}+a_{1}h+\cdots +a_{n}h^{n}$ is a
polynomial function such that $f(x+h)-P(h)=h^{n}\eta (h)$ where $\eta (h)$
tends to $0$ with $h$, then the \emph{Peano derivative of order }$n$ is $%
n!a_{n}$ (see, for example,  Weil \cite[(1995)]{Weil}, Svetic and
Volkmer \cite[(1998)]{svetic}). Peano gives an example of function that is
discontinuous in every neighborhood of $x$, and for which the Peano
derivatives of all order exist.},

\item derivation of measures in Peano \cite[(1887) p.\,169]{peano87} (see also Greco \cite%
{Greco}), and

\item mean value theorem.
\end{enumerate}

Here is an excerpt concerning the latter. In \emph{Calcolo Geometrico} \cite[%
(1888)]{calcolo88} Peano gives a \emph{mean value theorem} for vector-valued
functions $f$ of one variable \footnote{%
Peano gives a proof (by scalarization) in \cite[(1895) p.\,975]{saggio}.}, that is,
if $f$ has an $(n+1)$-derivative $f^{(n+1)}$ on $[t,t+h]$, then there exists
an element $k$ of the closed convex hull of the image of interval $[t,t+h]$
by $f^{(n+1)}$ such that%
\begin{equation}
f(t+h)=f(t)+hf^{\prime }(t)+\cdots +\dfrac{h^{n}}{n!}f^{(n)}(t)+\dfrac{%
h^{n+1}}{(n+1)!}k.
\end{equation}

Here is another surprise in front of the first appearance of modern notions
of \emph{convex set} \footnote{%
Peano employs the concept of convex set for the first time to axiomatic
foundation of geometry \cite[(1889) p.\,90 Axioma XVII]{peano_principii}; more
precisely, his \emph{Axiom XVII of continuity} states: \emph{Let $A$ be a
convex set of points, and let $x$, $y$ points such that $x\in A$ and $%
y\notin A$. Then there exists a point $w\in xy$ (the open segment between $x$
and $y$) such that $xw\subset A$ and $wy\cap A=\varnothing $}.} and \emph{%
convex hull}, as we thought that it was Minkowski who introduced these
concepts for the first time in \cite[(1896)]{minkowski}.

Among the first who studied the modern notion of differentiability of
functions of several variables, were  Stolz \cite[(1893), \S .IV\,8
p.\,130]{stolz},  Pierpont \cite[(1905) p.\,269]{pierpont}, W.H.
Young \cite[(1910) p.\,21]{young}. Besides, Maurice Fr\'{e}chet in \cite[%
(1911)]{Frechet_diff}\footnote{%
A month  after the publication of \cite[(1911)]{Frechet_diff} in which he
presented the concept of differentiability, Fr\'{e}chet publishes a second
Note \cite[(1911)]{Frechet_diff2} in order to recognize the priority of Young.}
defines the \emph{differential} (of a function of two variables) as the
linear part of the approximation%
\begin{equation}
f(a+h,b+k)-f(a,b)=hp+kq+h\rho (h,k)+k\sigma (h,k),
\end{equation}%
where $\rho $ and $\sigma $ tend to $0$ with $h$ and $k$.

Peano's definition liberates the concept of derivative from the coordinate
system and from partial derivatives. The definitions (\ref{differ})-(\ref%
{der+vect}) are coordinate-free in contrast with the predominant habit of
the epoch.

We have found no clear evidence in the mathematical literature of an
acknowledgement of Peano's definition of derivative, with the exception of a
paper \cite[(1921)]{wikosz} of Wilkosz, where Peano is cited jointly with
Stolz.

\section{Limits of variable sets}

In \emph{Applicazioni Geometriche} (p.\,30) Peano introduces a notion of
limit of straight lines, planes, circles and spheres (that depend on
parameter). He considers these objects as sets, which leads him to extend
the definition of limit to variable figure (in particular, curves and
surfaces).

A \emph{variable figure} (or \emph{set}) is a family, indexed by the reals,
of subsets $A_{\lambda }$ of an affine Euclidean space $X$. Peano defines in 
\emph{Applicazioni Geometriche} \cite[(1887) p.\,302]{peano87} the \emph{%
lower limit of a variable figure} by%
\begin{equation}
\Li_{\lambda \rightarrow +\infty }A_{\lambda }:=\{y\in
X:\lim_{\lambda \rightarrow +\infty }\dist(y,A_{\lambda
})=0\}.
\end{equation}

In the edition of \emph{Formulario Mathematico} \cite[(1908) p.\,237]%
{peano_1908} we find the lower limit together with a definition of \emph{%
upper limit of a variable figure}:%
\begin{equation}
\Ls_{\lambda \rightarrow +\infty }A_{\lambda }:=\{y\in
X:\liminf_{\lambda \rightarrow +\infty }\dist(y,A_{\lambda
})=0\},
\end{equation}%
which we have already seen in (\ref{upper+limit}). Besides, he writes down
(p.\,413) the upper limit as%
\begin{equation}
\Ls_{n\rightarrow \infty }A_{n}=\bigcap_{n\in \mathbb{N}}%
\cl\bigcup_{k\geq n}A_{k},
\end{equation}%
where the \emph{closure }$\cl A$ of a set $A$, is defined
(p.\,177) by\footnote{%
Peano defined closure, interior and boundary earlier in \cite[(1887)
pp.\,152-158]{peano87}; later, these notions were introduced by Jordan in 
\cite[(1893)]{Jordan2}. Peano relates the closure with the concept of \emph{%
closed set} of Cantor:\emph{\ the closure of }$A$\emph{\ is the least closed
set including }$A$.} 
\begin{equation}
\cl A:=\{y\in X:\dist (y,A)=0\}.
\end{equation}

In several papers, Peano analyzes the meanings that are given in mathematics
to the word \emph{limit} (see, for example, \cite[(1894)]{peano_Amer}):
least upper bound, greatest lower bound of a set, (usual) limit and adherence
of sequences and functions.

Peano conceives the ``upper limit of variable sets'' as a natural extension
of the \emph{adherence} of functions. He attributes to Cauchy the
introduction of adherence, see \cite[(1894) p.\,37]{peano_Amer} where he
says:

\begin{quotation}
Selon la d\'efinition de la limite, aujourd'hui adopt\'ee dans tous les
trait\'es, toute fonction a une limite seule, ou n'a pas de limite. [\dots]

Cette id\'{e}e plus g\'{e}n\'{e}rale de la limite [\emph{the adherence}] est
clairement \'{e}nonc\'{e}e par Cauchy; on lit en effet dans son Cours
d'Analyse alg\'{e}brique, 1821, p.\,13: \guillemotleft %
Quelquefois...une expression converge \`{a}-la-fois vers plusi\-eurs limites
diff\'{e}rentes les unes des autres\guillemotright , et \`{a} la page 14 il
trouve que les valeurs limites de $\sin \frac{1}{x}$, pour $x=0$,
constituent l'intervalle de $-1$ \`{a} $+1$. Les auteurs qui ont suivi
Cauchy, en cherchant de pr\'{e}ciser sa d\'{e}finition un peu vague, se sont
mis dans un cas particulier.
\end{quotation}

Peano studies the notion of \textquotedblleft lower limit of variable
sets\textquotedblright , in particular, in a celebrated article on existence
of solutions of a system of ordinary differential equations \cite[(1888)]%
{peano_equa_diff}. Peano carries on the proof of existence in a framework of
logical and set-theoretic ideography, thanks to which he is able to detect
the \emph{axiom of choice}\footnote{%
Peano proves the existence of a solution with the aid of approximated
solutions. In order to obtain a solution, he is confronted with a problem of
non-emptiness of the lower limit of a sequence of subsets of a finite-dimensional Euclidean
space. To this end, he needs to select an element from every set of the
sequence. At that point he realizes that he would need to make \emph{%
infinite arbitrary choices}, which, starting from the paper of Zermelo \cite%
{zermelo} of 1904 is called \emph{Axiom of choice}. He avoids to apply a new
axiom, which is not present in mathematical literature and, consequently,
the tradition does not grant it. Instead, using the lexicographic order, he
is able to construct a particular element of every set, because the sets of
the sequence are compact.}.

The awareness of the problem of \textquotedblleft limits of variable
sets\textquotedblright\ was present on the threshold of the 20-th Century
(for example see Manheim \cite[(1964)]{manheim}).\footnote{%
In \cite[(1903)]{Borel} (see also Manheim \cite[(1964) p.\,114]{manheim}) 
Borel suggests a \textquotedblleft promising\textquotedblright\ notion of
limit of straight lines and of planes, that is, 16 years after the
introduction of the limit of arbitrary sets in \emph{Applicazioni Geometriche}%
.} The book of Kuratowski \cite[(1948) p.\,241]{kuratowski} \footnote{%
Kuratowski, by his work, consacrates the use of upper and lower limits in
mathematics, that are called today \emph{upper and lower Kuratowski limits.}}
has definitely propagated the concept of limit of variable sets.

Among the first mathematicians who studied the limits of variable sets are
 Burali-Forti \cite[(1895)]{burali-forti} \footnote{%
Burali-Forti studies only lower limits.}, Zoretti \cite[(1905) p.\,8]%
{zoretti_J} and \cite[(1909)]{zoretti_B}, 
  Janiszewski \cite[(1911)]{janiszewski}, 
  Hausdorff \cite[(1914) p.\,234]{hausdorff} \footnote{%
In the celebrated {\em Grundz{\" u}ge der
Mengelehre} (1914) Hausdorff studies both upper and lower limits. Moreover, he defines a metric
on the set of bounded subsets of a metric space $X$ (\emph{Hausdorff distance}%
) and proves that the related convergence of bounded subsets $\{A_{\lambda}\}_{\lambda}$   to $A$ (as $\lambda\to+\infty$) is equivalent to $\Ls_{\lambda
\rightarrow +\infty }A_{\lambda }\subset A\subset \Li_{\lambda
\rightarrow +\infty }A_{\lambda }$, if $X$ is compact. 
}, Vietoris \cite[1922]{vietoris}\footnote{Kuratowski limits, Hausdorff distance and Vietoris's topology (see Reitberger \cite[(2002) p.\,1234]{reitberger}) are milestones  in the search of notions of limit of variable sets.}, 
Vasilesco \cite[(1925)]{vasilesco},  Cassina \cite[(1926-27)]{Cassina} and Kuratowski \cite[(1928)]{kuratowski28}.

\section{Tangent cones}

In \emph{Applicazioni Geometriche} (pp.\,58, 116) Peano gives a
metric definition of tangent straight line and tangent plane, then reaches,
in a natural way, a unifying notion: that of \emph{affine tangent cone}:%
\begin{equation}
\tang(A,x):=\Li_{\lambda \rightarrow +\infty }\left(
x+\lambda (A-x)\right) .
\end{equation}%
Later, in \emph{Formulario Mathematico} (p.\,331), he introduces another type
of tangent cone, namely%
\begin{equation}
\Tang(A,x):=\Ls_{\lambda \rightarrow +\infty }\left( x+\lambda
(A-x)\right)
\end{equation}%
To distinguish the two notions above, we shall call the first \emph{lower
affine tangent cone} and the second \emph{upper affine tangent cone} 
\footnote{%
Of course, $\tang(A,x)$ is defined with the aid of lower limit, while $%
\Tang(A,x)$ with the aid of the upper limit of the same homothetic
sets. Hence, $\tang(A,x)\subset \Tang(A,x)$. For the
covenience of the reader, in order to compare the two definitions, we give
their alternative descriptions in terms of limits of sequences:%
\begin{eqnarray}
\qquad\quad\tang(A,x) =x+\left\{ v:\exists \left\{ x_{n}\right\} _{n}\subset A%
\text{ such that }x=\lim_{n\to\infty} x_{n}\text{ and }v=\lim_{n\to\infty} \dfrac{x_{n}-x}{1/n}%
\right\} ;\label{tan+new}
\end{eqnarray}
\begin{eqnarray*}
\Tang(A,x) =x+\left\{ v:\exists \left\{ \lambda_{n}\right\}
_{n}\rightarrow 0^{+},\;\exists \left\{ x_{n}\right\} _{n}\subset A\text{
such that }x=\lim_{n\to\infty} x_{n}\text{ and }v=\lim_{n\to\infty} \dfrac{x_{n}-x}{\lambda _{n}}%
\right\} .
\end{eqnarray*}%
The second formula is standard, while we have never seen in the literature
the first one (\ref{tan+new}). We have not found in Peano's papers any example of set $A$
for which the cones above are different. Here is another, perhaps most
intuitive, formula for the "lower" affine tangent cone:%
\begin{equation}
\tang(A,x)=x+\left\{ v:\exists \gamma :[0,1]\rightarrow A\text{ such
that }x=\gamma (0),\gamma ^{\prime }(0)\text{ exists and }v=\gamma ^{\prime
}(0)\right\} .
\end{equation}%
\par
Notice that $\tang(A,x)=\Tang(A,x)$ in case of differential
manifold $A$ (at $x$).}.\emph{\ }Peano lists several properties of the upper
tangent cone.
If $A$ is a subset of a Euclidean affine space $X$%
, then

\begin{enumerate}
\item If $x\notin \cl A$ then $\Tang(A,x)=\varnothing ;$

\item If $x$ is isolated in $A$ then $\Tang(A,x)=\left\{ x\right\}
; $

\item If $x\in \cl (A\backslash \left\{ x\right\} )$ then $
\Tang(A,x)\neq \varnothing ;$

\item If $x\in \intt A$ then $\Tang(A,x)=X;$

\item If $y\in \Tang(A,x)\setminus\{x\}$ then $x+\mathbb{R}%
_{+}(y-x)\subset \Tang(A,x);$

\item If $A\subset B$ then $\Tang(A,x)$ $\subset \Tang%
(B,x);$

\item $\Tang(A\cup B,x)=\Tang(A,x)\cup \Tang%
(B,x); $

\item $\Tang(\Tang(A,x),x)=\Tang(A,x)$.\footnote{%
One find in fourth edition of \emph{Formulaire Math\'{e}matique} \cite[%
(1903) p.\,296]{peano_F4} the same definition of upper affine tangent cone
and, besides, the eight properties (1)--(8). Besides, one find both
lower and upper limit of variable sets \cite[(1903) p.\,289]{peano_F4}.

We have not found in any of five editions of \emph{Formulario Mathematico} (= collection of logical and set-teoretical formulas) any other property on  tangent cones. Today other fundamental properties are well-known: (1) $x\in \cl  A \iff \Tang(A,x)\ne\emptyset\iff x\in \Tang(A,x)$; (2) $x\in \cl (A\setminus\{x\})\iff  \Tang(A,x)\setminus \{x\}\ne\emptyset$; $(3)$ $\Tang(A,x)=\Tang(A\cap B,x)$, if $x\in \intt B$; $(4)$  $\Tang(A,x)=\Tang\big(\cl (A),x\big)$; finally, $(5)$ $\Tang(A,x)$ is closed (because it is an upper limit of variable sets).}
\end{enumerate}

As usual, after abstract investigation of a notion, Peano considers
significant special cases; he calculates the upper affine tangent cone in
several basic figures (closed ball, curves and surfaces parametrized in a
regular way).

Various types of tangent cones have been studied in the literature. Their
definitions depend on variants of limiting process. The most known
contribution to the investigation of tangent cones is due to 
Bouligand \cite[(1932) p.\,60]{bouligand}. One can find a mention about other
contributors in a paper \cite[(1937) p.\,241]{Frechet_3} of Fr\'{e}chet:

\begin{quotation}
Cette th\'{e}orie des \textquotedblleft contingents et
paratingents\textquotedblright\ dont l'utilit\'{e} a \'{e}t\'{e} signal\'{e}%
e d'abord par M. Beppo Levi, puis par M. Severi, mais dont on doit \`{a} M.
Bouligand et \`{a} ses \'{e}l\`{e}ves d'avoir entrepris l'\'{e}tude syst\'{e}%
matique.
\end{quotation}

The diffusion of the concept of tangent cone was due mainly to 
Saks \cite[(1937) pp.\,262--263]{Saks}, who adopted the definition of
Bouligand, and to Federer \cite[(1959) p.\,433]{Federer}, who introduced it
in a modern vector version: if $x\in A$ then define the \emph{upper vector tangent 
cone}: 
\begin{equation}
\left. 
\begin{array}{c}
\Tan(A,x):=\left\{ 0\right\} \cup \\ 
\left\{ u\neq 0:\forall \varepsilon >0,\;\exists y\in
A,\;0<|y-x|<\varepsilon ~\text{and }\left| \dfrac{y-x}{| y-x| }-%
\dfrac{u}{| u| }\right| <\varepsilon \right\} .%
\end{array}%
\right.  \label{tan+fed}
\end{equation}%
Federer does not give any reference of the origin of the definition (\ref%
{tan+fed})\footnote{%
The book of Saks \cite[(1937)]{Saks} is among bibliographic references in Federer \cite[(1959)]{Federer}.}. Notice
that 
\begin{equation}
\Tang(A,x)=x+\Tan(A,x).
\end{equation}

Neither  Whitney cites nobody in \cite[(1972) chap.\,7]%
{whitney}, where he introduces six variants of vector tangent cone among which, one
recognizes the upper vector tangent cone discussed above (\ref{tan+fed}).

We can say that, as far as tangent cones are concerned, main references are,
respectively, Bouligand in optimization theory, Ferderer in geometric
measure theory and calculus of variations and Whitney in differential
geometry. A rare direct reference to Peano's definition is that of Guido
Ascoli \footnote{%
One should not confound Guido Ascoli(1887-1957) with Giulio
Ascoli(1843-1896), the latter known because of the Ascoli-Arzel\`{a}
theorem. Guido proved the geometric version \cite[(1933)]{ascoliB} of the
Hahn-Banach theorem for separable normed spaces; a year later,
Mazur proved it for arbitrary normed spaces \cite[(1933)]{mazurB}.} \cite[(1953$)$]%
{ascoli}, who writes about Peano's work in \cite[(1955) pp.\,26-27]%
{Ascoli_inbook}:

\begin{quotation}
[...] il merito maggiore [...] specialmente delle \emph{Applicazioni} [\emph{%
Geometriche}], non sta tanto nel metodo usato, quanto nel contenuto; ch\'{e}
vi sono profusi, in forma cos\`{\i} semplice da parere definitiva, idee e
risultati divenuti poi classici, come quelli sulla misura degli insiemi,
sulla rettificazione delle curve, sulla definizione dell'area di una
superfice, sull'integrazione di campo, sulle funzioni additive d'insieme; ed
altri che sono tuttora poco noti o poco studiati. Ci basti indicare tra
questi il concetto di limite di una figura variabile, destinato a
ricomparire, con altro nome di autore, quarant'anni dopo presso la scuola di
"geometria infinitesimale diretta" del Bouligand, e l'originalissima
definizione di "figura tangente ad un insieme in un punto", che ha fornito a
chi scrive, or \`{e} qualche anno, la chiave di una difficile questione
asintotica.
\end{quotation}

The \emph{contingent cone} of Bouligand \cite[(1932)]{bouligand} is defined
by Saks in \cite[(1937) p.\,262]{Saks} as follows: if $x$ is an accumulation
point of $A$, then the \emph{contingent cone} of $A$ at $x$ is given by%
\begin{equation}
\Cont(A,x):=\left\{ l:l\text{ tangent half-line to }A\text{ at }%
x\right\}
\end{equation}%
where a \emph{half-line} $l$ issued from $x$ is said to be \emph{tangent to }%
$A$ at $x$ if{} there exist a sequence $\left\{ x_{n}\right\} _{n}\subset $ $%
A$ and a sequence of half-lines $\left\{ l_{n}\right\} _{n}$ issued from $x$
such that $x\neq x_{n}\in l_{n},x=\lim x_{n}$ and the angle between $l_{n}$
and $l$ tends to $0$.

Peano's upper affine tangent cone, Federer's upper vector tangent cone and
Bouligand's contingent cone describe the same intuitive concept in terms,
respectively, of points of affine space (via \emph{blow-up}), of vectors
(via \emph{directions} of tangent half-lines) and of half-lines (via \emph{%
limits} of half-lines)\footnote{Bouligand, in spite of his knowledge of vector
spaces (see, for example, \cite[(1924)]{Boul_Vect} and his introduction to
the French translation of \cite[(1918)]{Weyl} of  Weyl, does not use
vectors while defining the contingent cone. In the preface to \cite[1957]{felix} Bouligand appraises Peano's \emph{Calcolo Geometrico} \cite[(1888)]{calcolo88}: \guillemotleft Pour \^etre moins incomplet, il faudrait encore citer l'expos\'e repris par Peano en 1886 [sic!]  \emph{du calcul extensif} de Grassmann, [et] l'article fondamental malgr\'e sa bri\`evet\'e paru en 1900 [sic!]  dans l'Enseignement Math\'ematique au sujet des relations d'\'equivalence, r\'edig\'e par Burali-Forti (Sur quelques notions d\'eriv\'ees de la notion d'\'egalit\'e et leurs applications dans la science).\guillemotright\ Burali-Forti was Peano's assistant and friend;  the article quoted by Bouligand is \cite[(1899)]{burali-fortiE}.}. Finally, observe that the tangent cone is built on the notion of distance  by Peano, of norm by Federer and of angle (consequently, of scalar product) by Bouligand. 

\section{Maxima and minima}

In \emph{Applicazioni Geometriche} (pp.\,143-144) Peano analyzes the
variation of a real-valued function around a point in a particular direction%
\textbf{\ }$\mathbf{p}$ in terms of the scalar product of the derivative at
that point with $\mathbf{p}$.

\begin{theorem}
\label{optimal 87} Let $f$ be a real-valued function such that $Df(\bar{x}%
)\neq 0$. Let $\mathbf{p}$ be a unit vector and $\left\{ x_{n}\right\} _{n}$
be a sequence  so that%
\begin{eqnarray} 
 &\bar x=\lim\nolimits_{n\to\infty} x_{n}\text{\qquad and\qquad}
\mathbf{p}=\lim\nolimits_{n\to\infty} \dfrac{x_{n}-\bar{x}}{\Vert x_{n}-\bar{%
x}\Vert }.  \label{p} \\
\text{If }&\langle Df(\bar{x}),\mathbf{p}\rangle >0\text{, then  }f(x_{n})>f(x_{0}) \text{ for almost all }n;\\
\text{if }&\langle Df(\bar{x}),\mathbf{p}\rangle <0\text{ then } f(x_{n})<f(x_{0})
\text{ for almost all }n. \quad\Box
\end{eqnarray}

\end{theorem}

Peano specifies that $\left\{ x_{n}\right\} _{n}$ in Theorem \ref{optimal 87}
can be taken either arbitrarily or constrained by some conditions, for
example, lying on a line or on a surface.

If $\left\{ x_{n}\right\} _{n}$ is included in a set $A$, then $\mathbf{p}$
is one of the directions (unitary vectors) of the upper vector tangent cone of $A$ at $x$%
. By taking all such sequences, we get all the directions of the upper vector
tangent cone (\ref{tan+fed}) of $A$ at $x$. Hence, by relating (\ref{p}) to
upper vector tangent cone, Theorem \ref{optimal 87} implies\smallskip

\begin{theorem}[Regula \textbf of Maximality] {If }$f${\ is differentiable
at }$x\in A${\ and }$f(x)=\max \{f(y):y\in A\}${, then}$\langle
Df(x),y-x\rangle \leq 0{\ for\ every\ }y\in \Tang(A,x).$
\end{theorem}

\smallskip\begin{theorem}[Regula \textbf of Minimality] {If }$f$ {\ is differentiable
at }$x\in A${\ and }$f(x)=\min \{f(y):y\in A\}${, then}$\langle
Df(x),y-x\rangle \geq 0{\ for\ every\ }y\in \Tang(A,x)$
\end{theorem}
\smallskip

One finds both Theorems 6.2 and 6.3 in \emph{Formulario Mathematico} (p.\,335). It is worthwhile to note that Peano's use of \emph{Regula} exhibits the
normality of gradient with respect to the constraint.

Optimization problems were among principal interests of Peano. His research
with regard to these problems was intense, continual and influential. The
precision with which Peano studied maxima and minima was notorious.

Hancock, student of Weierstrass, is  author of a booklet: \emph {Lectures on the theory of maxima and minima of functions of several variables.\,Weierstrass' theory} (1903).  In the second edition of this book  he says \cite[(1917) pp.\,iv-v]{Hancock}:

\begin{quotation}
In the preface to the German translation by Bohlmann and Schepp of Peano's
of \emph{Calcolo differenziale e principii di calcolo integrale}, Professor
A. Mayer [editor of Math.\thinspace Annalen together with Felix Klein]
writes that this book of Peano not only is a model of precise presentation
and rigorous deduction, whose propitious influence has been unmistakably
felt upon almost every calculus that has appeared (in Germany) since that
time (1884), but by calling attention to old and deeply rooted errors, it
has given an impulse to new and fruitful development.

The important objection contained in this book [\emph{Calcolo differenziale
e principii di calcolo integrale}] (Nos. 133-136) showed unquestionably that
the entire former theory of maxima and minima needed a thorough renovation;
and in the main Peano's book is the original source of the beautiful and to a great
degree fundamental works of Scheeffer, Stolz, Victor v.\,Dantscher, and
others, who have developed new and strenuous theories for \emph{extreme}
values of functions. Speaking for the Germans, Professor A. Mayer, in the
introduction to the above-mentioned book, declares that there has been a
long-felt need of a work which, for the first time, not only is free from
mistakes and inaccuracies that have been so long in vogue but which,
besides, so incisively penetrates an important field that hitherto has been
considered quite elementary.
\end{quotation}

\section{Appendix}

All articles of Peano are collected in \emph{Opera Omnia} \cite{peano_omnia}, a compact disk read only memory (CD-ROM),  edited by S. Roero. Selected works of Peano were assembled and commented in \emph{Opere scelte} \cite[(1957-59)]{peano_opere} by Cassina, a student of Peano. A few have English translations in \emph{Selected Works} \cite[(1973)]{peano_english}. Regrettably, even fewer Peano's articles have a public URL and are freely downloadable. 

One finds the following articles of Peano

 \bigskip
 
 \noindent\begin{tabular}{ll}
 \hline
 in \emph{Opere scelte}, vol. 1:& \begin{minipage}{180pt}\cite[(1884)]{peano_jordan}, \cite[(1884)]{peano_gilbert}, \cite[(1888)]{peano_ED},
 
  \cite[(1890)]{peano_area}, \cite[(1890)]{peano_equa_diff}, \cite[(1891)]{peano_taylor}, 
 
 \cite[(1892)]{peano_mathesis}, \cite[(1894)]{peano_exp},  \cite[(1894)]{peano_Amer},
 \end{minipage}\cr
\hline

in \emph{Opere scelte}, vol. 2: &\cite[(1889)]{peano_principii}, \cr
\hline

in \emph{Opere scelte}, vol. 3: & \cite[(1896)]{saggio}, \cite[(1898)]{peano_vettori}, \cr
\hline

 in \emph{Selected Works}: & \begin{minipage}{180pt}\cite[(1887) pp.\,152--160, 185--7]{peano87},
 
  \cite[(1887) pp.\,1--32]{calcolo88}, \cite[(1890)]{peano_area}, 
  
  \cite[(1894)]{peano_exp}, 
 \cite[(1896)]{saggio}. \end{minipage}\cr
\hline
 
 \end{tabular}

\bigskip
Before the bibliography we attach several pages of 
 \emph{Applicazioni
Geometriche} (AG) and of \emph{Formulario Mathematico} (FM) corresponding to 
Regula, limits of variable sets, derivatives and tangent cones.

For reader's convenience, we provide a chronological list of some
mathematicians mentioned in the paper, together with biographical sources.

Every \texttt{html} file listed below can be attained at University of St
Andrews's web-page \url{http://www-history.mcs.st-and.ac.uk/history/}
\medskip

\textsc{Jacobi}, Carl (1804-1851), see 
\href{http://www-history.mcs.st-and.ac.uk/history/Biographies/Jacobi.html}{St
Andrews's web-page}

\textsc{Hamilton}, William R. (1805-1865), see 
\href{http://www-history.mcs.st-and.ac.uk/history/Biographies/Hamilton.html}{St
Andrews's web-page}

\textsc{Grassmann}, Hermann (1809-1877), see 
\href{http://www-history.mcs.st-and.ac.uk/history/Biographies/Grassmann.html}{St
Andrews's web-page}

\textsc{Weierstrass}, Karl (1815-1897), see 
\href{http://www-history.mcs.st-and.ac.uk/history/Biographies/Weierstrass.html}{St
Andrews's web-page}

\textsc{Genocchi}, Angelo (1817-1889), see 
\href{http://www-history.mcs.st-and.ac.uk/history/Biographies/Genocchi.html}{St
Andrews's web-page}

\textsc{Jordan}, Camille (1838-1922), see 
\href{http://www-history.mcs.st-and.ac.uk/history/Biographies/Jordan.html}{St
Andrews's web-page}

\textsc{Mayer}, Adolph (1839-1908), see 
\href{http://www-history.mcs.st-and.ac.uk/history/Biographies/Mayer.html}{St
Andrews's web-page}

\textsc{Darboux}, Gaston (1842-1917), see 
\href{http://www-history.mcs.st-and.ac.uk/history/Biographies/Darboux.html}{St
Andrews's web-page}

\textsc{Stolz}, Otto (1842-1905), see 
\href{http://www-history.mcs.st-and.ac.uk/history/Biographies/Stolz.html}{St
Andrews's web-page}

\textsc{Ascoli}, Giulio (1843-1896), see May \cite[(1973) p.\,62]{may}

\textsc{Cantor}, Georg (1845-1918), see 
\href{http://www-history.mcs.st-and.ac.uk/history/Biographies/Cantor.html}{St
Andrews's web-page}

\textsc{Dini}, Ulisse (1845-1918), see 
\href{http://www-history.mcs.st-and.ac.uk/history/Biographies/Dini.html}{St
Andrews's web-page}

\textsc{Dantscher von Kollesberg}, Victor (1847-1921), see A.M.\,Monthly, 29
(1922)

\textsc{Klein}, Felix (1849-1925), see 
\href{http://www-history.mcs.st-and.ac.uk/history/Biographies/Klein.html}{St
Andrews's web-page}

\textsc{Poincar\'e}, Henri (1854-1912), see 
\href{http://www-history.mcs.st-and.ac.uk/history/Biographies/Poincare.html}{St
Andrews's web-page}

\textsc{Goursat}, Edouard (1858-1936), see 
\href{http://www-history.mcs.st-and.ac.uk/history/Biographies/Goursat.html}{St
Andrews's web-page}

\textsc{Peano}, Giuseppe (1858-1932), see  \cite{peano_vita} and  
\href{http://www-history.mcs.st-and.ac.uk/history/Biographies/Peano.html}{St
Andrews's web-page}

\textsc{Scheeffer}, Ludwig (1859-1885), see  Math. Annalen 26
(1886),p.\,197

\textsc{Burali-Forti}, Cesare (1861-1931), see 
\href{http://www-history.mcs.st-and.ac.uk/history/Biographies/Burali-Forti.html}{St
Andrews's web-page}

\textsc{Young}, William H. (1863-1942), see 
\href{http://www-history.mcs.st-and.ac.uk/history/Biographies/Young.html}{St
Andrews's web-page}

\textsc{Painlev\'e}, Paul (1863-1933), see 
\href{http://www-history.mcs.st-and.ac.uk/history/Biographies/Painleve.html}{St
Andrews's web-page}

\textsc{Pierpont}, James (1866-1938), see \emph{Bull.\,A.M.S.}\,45 (1939), p.\,481

\textsc{hancock}, Harris (1867-1944), see May \cite[(1973) p.\,185]{may}

\textsc{Hausdorf}, Felix (1868-1942), see 
\href{http://www-history.mcs.st-and.ac.uk/history/Biographies/Hausdorf.html}{St
Andrews's web-page}

\textsc{Borel}, Emile (1871-1956), see 
\href{http://www-history.mcs.st-and.ac.uk/history/Biographies/Borel.html}{St
Andrews's web-page}

\textsc{Carath\'{e}odory}, Constantin (1873-1950), see 
\href{http://www-history.mcs.st-and.ac.uk/history/Biographies/Caratheodory.html}{St
Andrews's web-page}

\textsc{Levi}, Beppo (1875-1961), see May \cite[(1973) p.\,238]{may}

\textsc{Fr\'{e}chet}, Maurice (1878-1973), see 
\href{http://www-history.mcs.st-and.ac.uk/history/Biographies/Frechet.html}{St
Andrews's web-page}

\textsc{Severi}, Francesco (1879-1961), see 
\href{http://www-history.mcs.st-and.ac.uk/history/Biographies/Severi.html}{St
Andrews's web-page}

\textsc{Zoretti}, Ludovic (1880-1948), see \texttt{http://catalogue.bnf.fr}

\textsc{Weyl}, Hermann (1885-1955), see 
\href{http://www-history.mcs.st-and.ac.uk/history/Biographies/Weyl.html}{St
Andrews's web-page}
 
\textsc{Ascoli}, Guido (1887-1957), see May \cite[(1973) p.\,63]{may}

\textsc{Janiszewski}, Zygmunt (1888-1920), see 
\href{http://www-history.mcs.st-and.ac.uk/history/Biographies/Janiszewski.html}{St
Andrews's web-page}

\textsc{Bouligand}, Georges (1889-1979), see \texttt{http://catalogue.bnf.fr}

\textsc{Vietoris}, Leopold (1891-2002), see Reitberger \cite{reitberger}

\textsc{Wilkosz}, Wiltold (1891-1941), see 
\href{http://www.wiw.pl/matematyka/Biogramy/}{Matematycy polscy}

\textsc{Banach}, Stefan (1892-1945), see 
\href{http://www-history.mcs.st-and.ac.uk/history/Biographies/Banach.html}{St
Andrews's web-page}

\textsc{Kuratowski}, Kazimierz (1896-1980), see 
\href{http://www-history.mcs.st-and.ac.uk/history/Biographies/Kuratowski.html}{St
Andrews's web-page}

\textsc{Cassina}, Ugo (1897-1964), see May \cite[(1973) p.\,222]{may}

\textsc{Vasilesco}, Florin (1897-1958),  see May \cite[(1973) p.\,368]{may}

\textsc{Saks}, Stanislaw (1897-1942), see 
\href{http://www-history.mcs.st-and.ac.uk/history/Biographies/Saks.html}{St
Andrews's web-page}

\textsc{Mazur}, Stanislaw (1905-1981), see 
\href{http://www-history.mcs.st-and.ac.uk/history/Biographies/Mazur.html}{St
Andrews's web-page}

\textsc{Whitney}, Hassler (1907-1989), see 
\href{http://www-history.mcs.st-and.ac.uk/history/Biographies/Whitney.html}{St
Andrews's web-page}

\textsc{Birkhoff}, Garrett (1911-1996), see 
\href{http://www-history.mcs.st-and.ac.uk/history/Biographies/Birkhoff_Garrett.html}{St
Andrews's web-page}

\begin{figure}[htbp]
\begin{center}
\includegraphics[scale=.5]{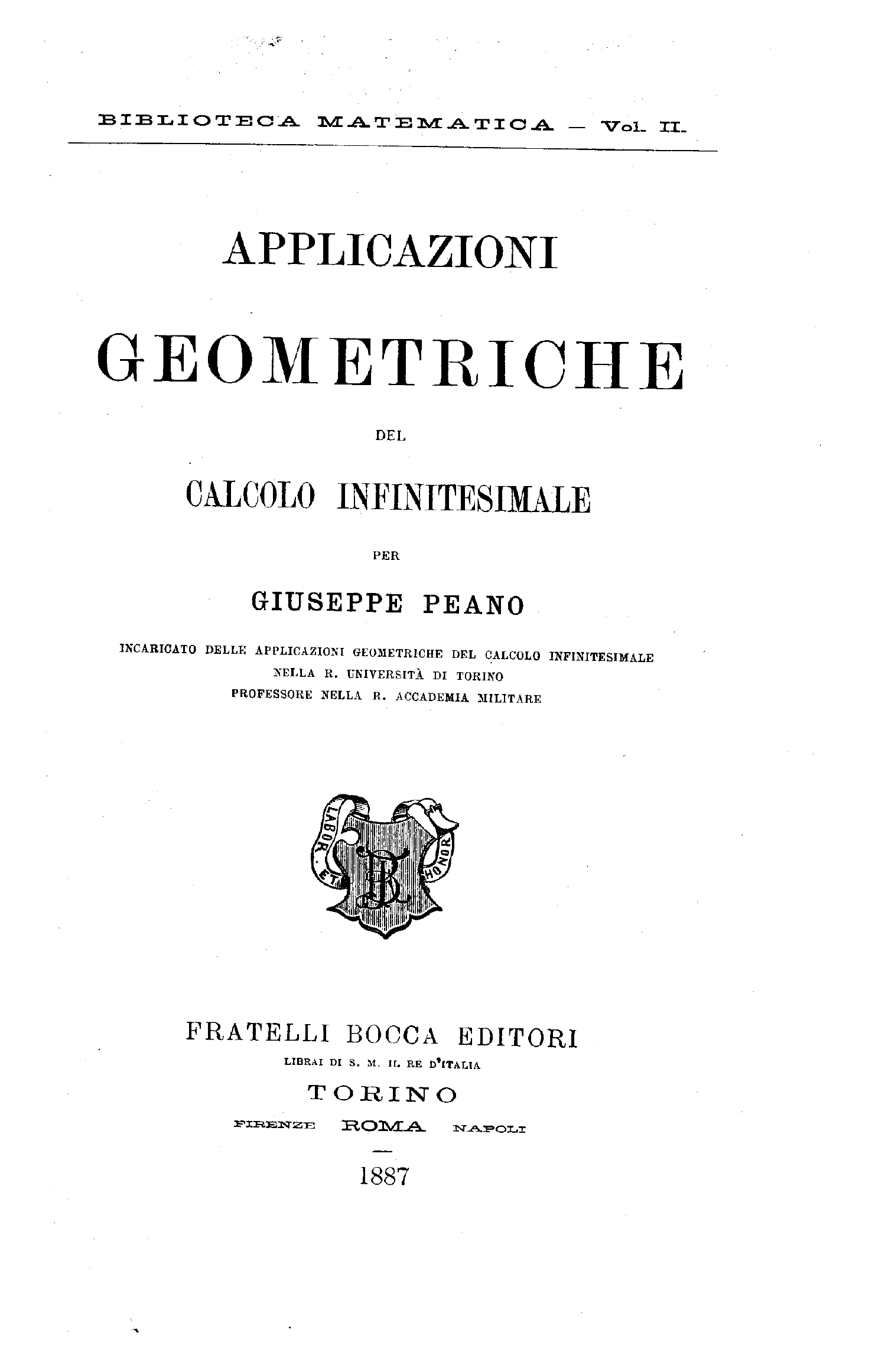}

\caption{Applicazioni geometriche}
\label{AGcopertina}
\end{center}
\end{figure}

\begin{figure}[htbp]
\begin{center}
\includegraphics[scale=.9]{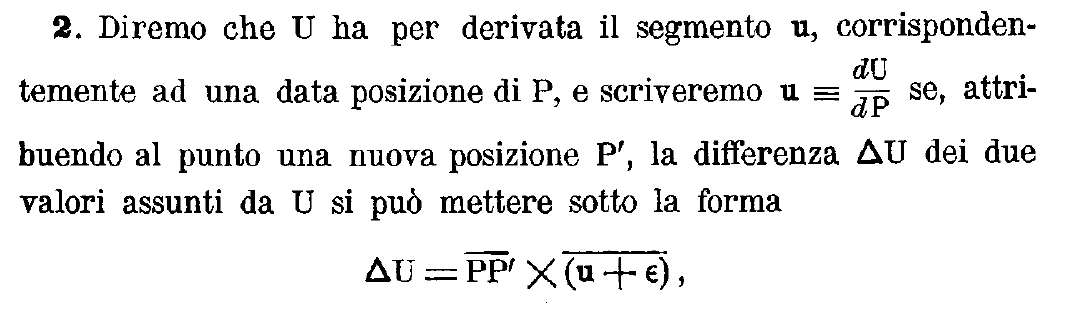}

\includegraphics[scale=.9]{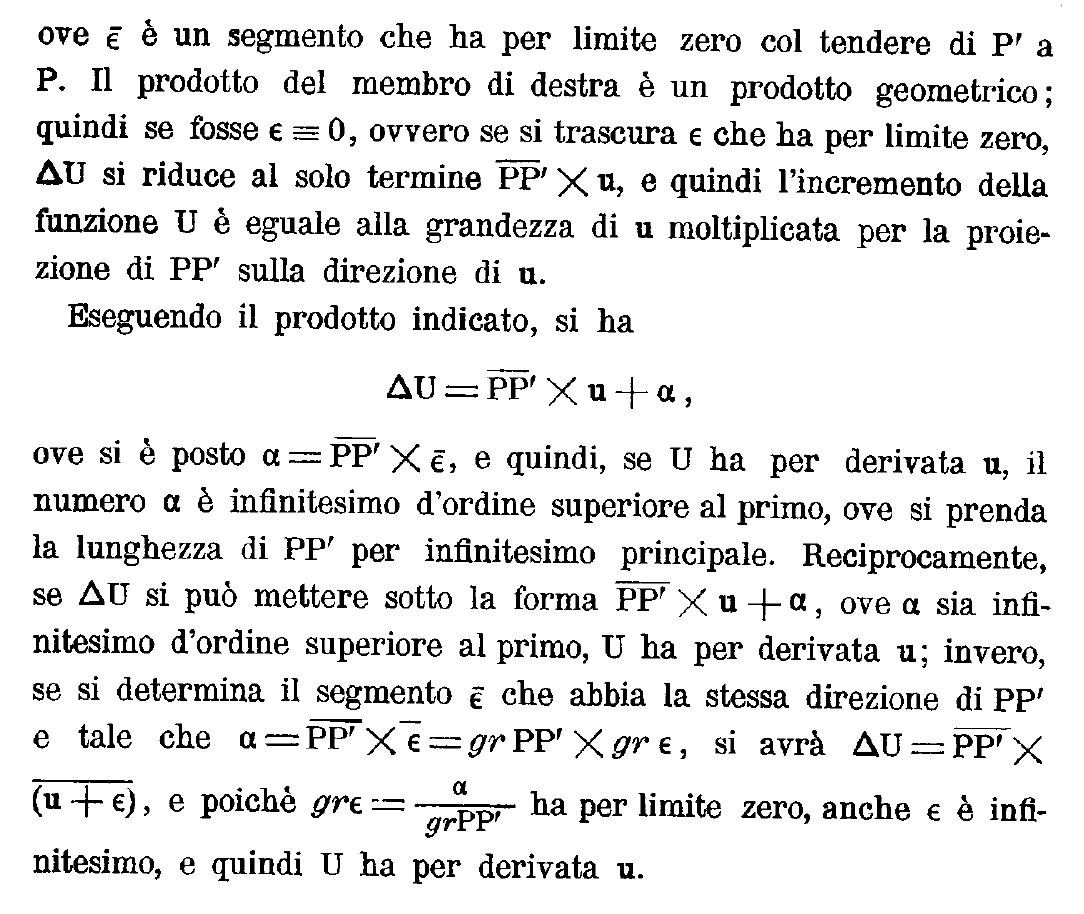}
\caption{AG (pp.\,131--132): derivative }
\label{AG131}
\end{center}
\end{figure}

\begin{figure}[htbp]
\begin{center}
\includegraphics[scale=.9]{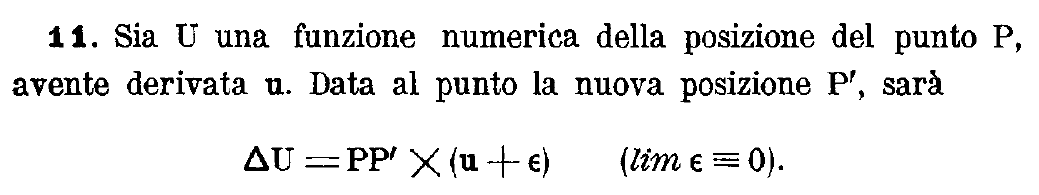}

\includegraphics[scale=.9]{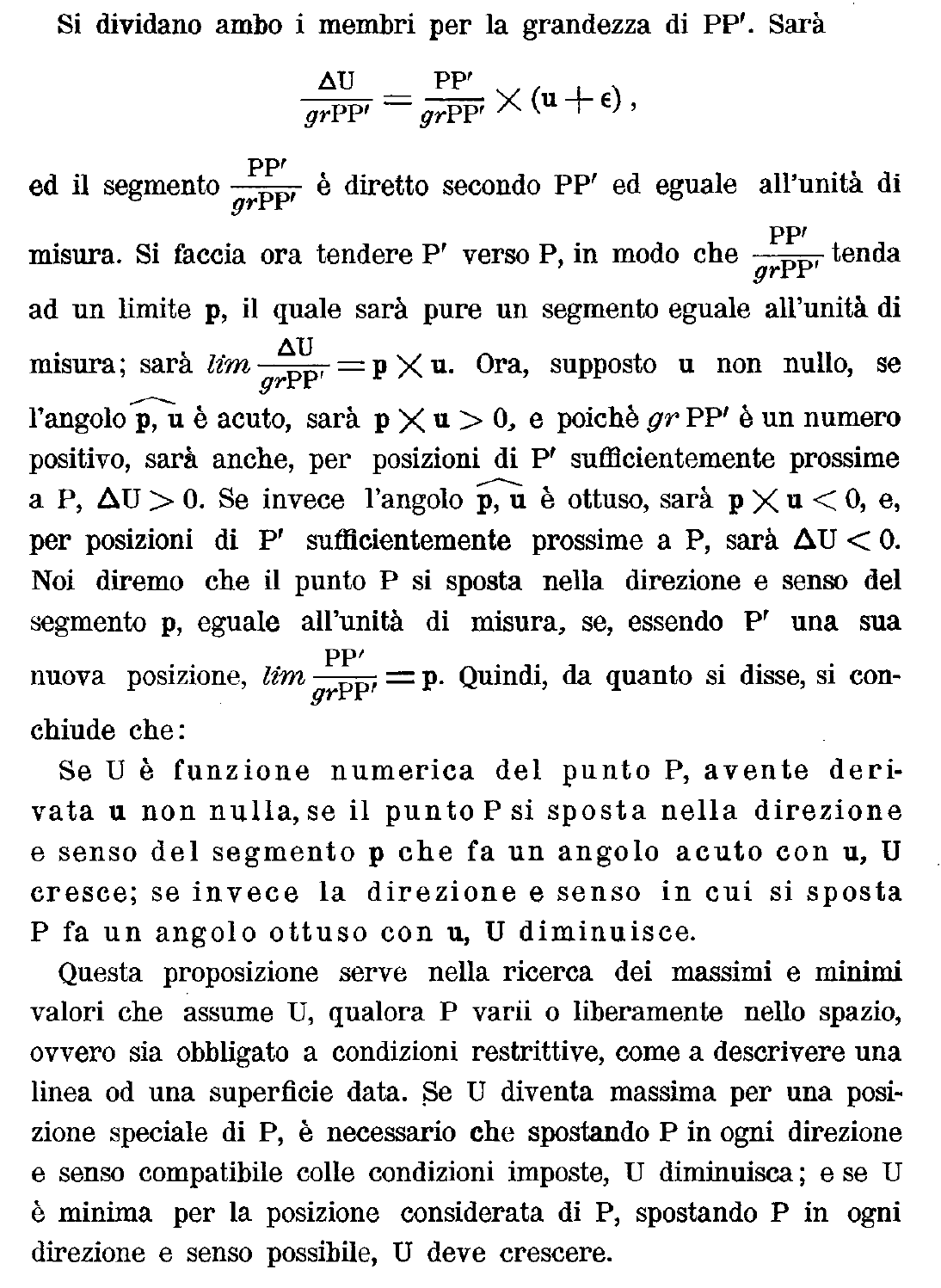}

\caption{AG (pp.\,143--144): Regula}
\label{AG143}
\end{center}
\end{figure}

\begin{figure}[htbp]
\begin{center}
\includegraphics[scale=.9]{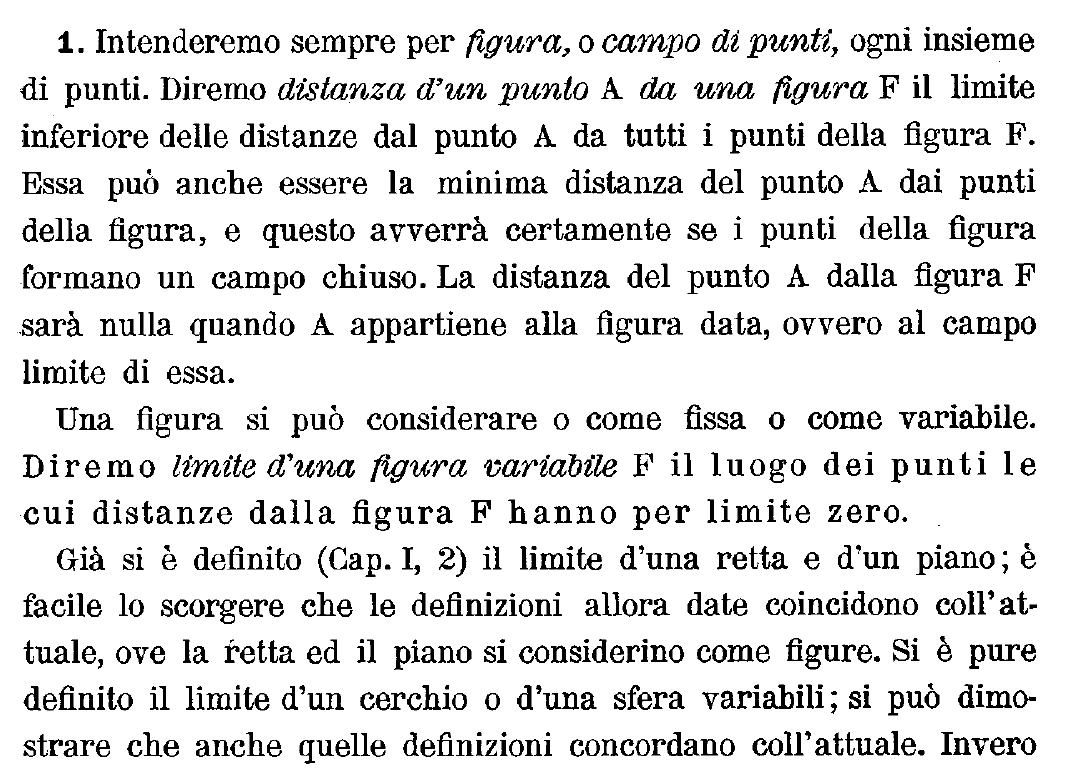}
\caption{AG (p.\,302): lower limit of variable sets}
\label{AG302}
\end{center}
\end{figure}

\begin{figure}[htbp]
\begin{center}
\includegraphics[scale=.9]{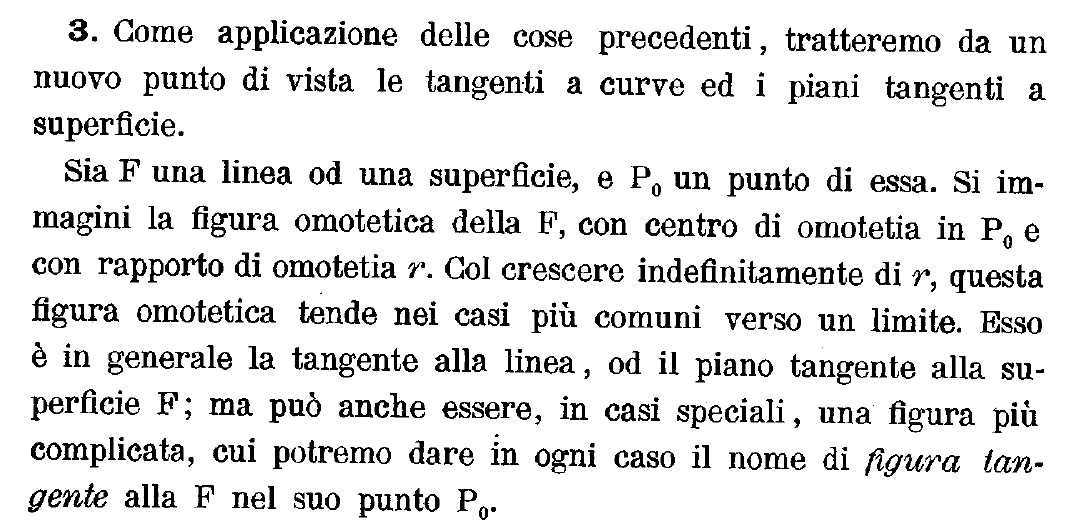}
\caption{AG (p.\,305): lower affine tangent cone}
\label{AG305}
\end{center}
\end{figure}

\begin{figure}[htbp]
\begin{center}
\includegraphics[scale=.5]{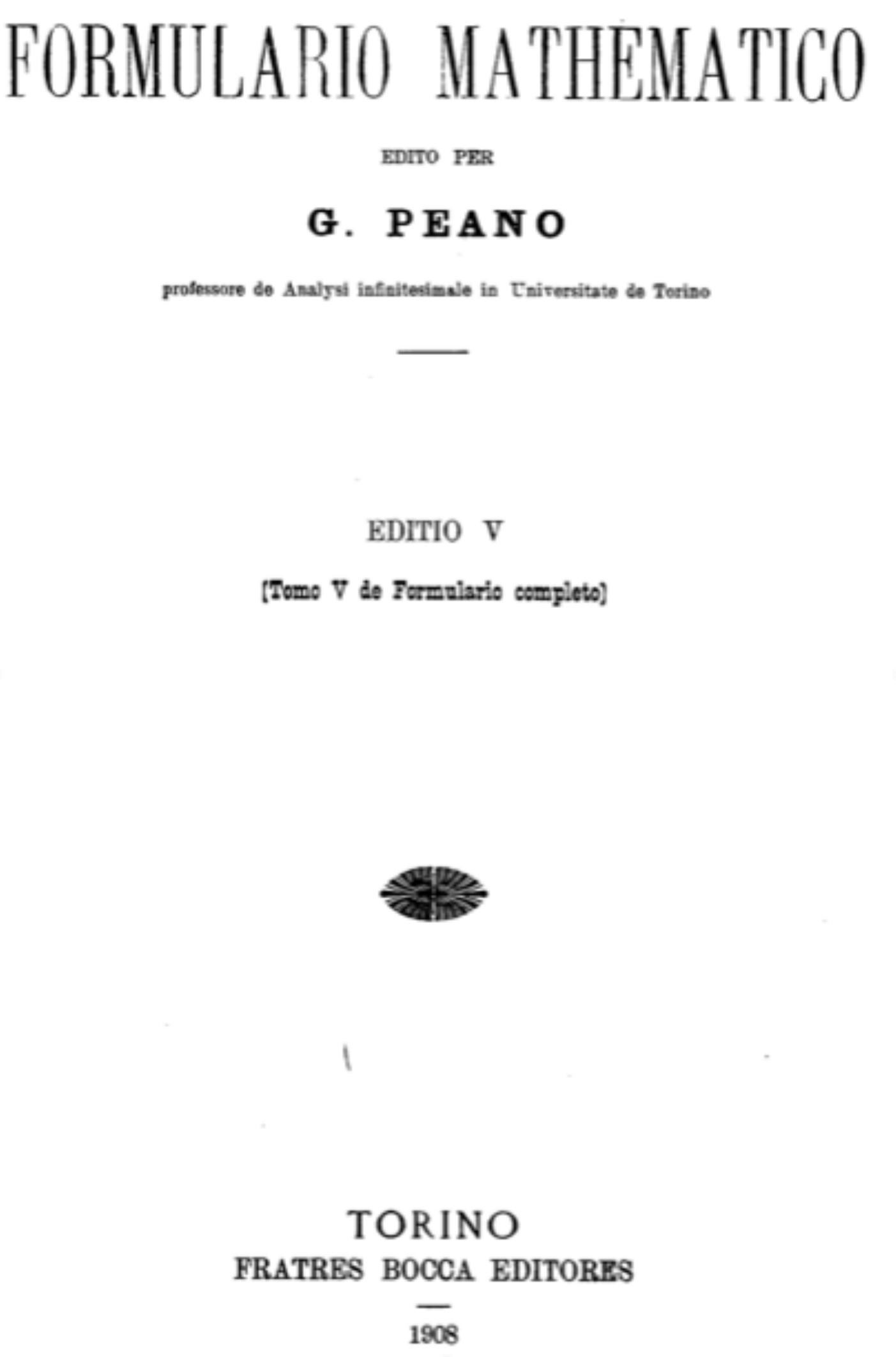}
\caption{Formulario Mathematico}
\label{FMcopertina}
\end{center}
\end{figure}

\begin{figure}[htbp]
\begin{center}
\includegraphics[scale=.83]{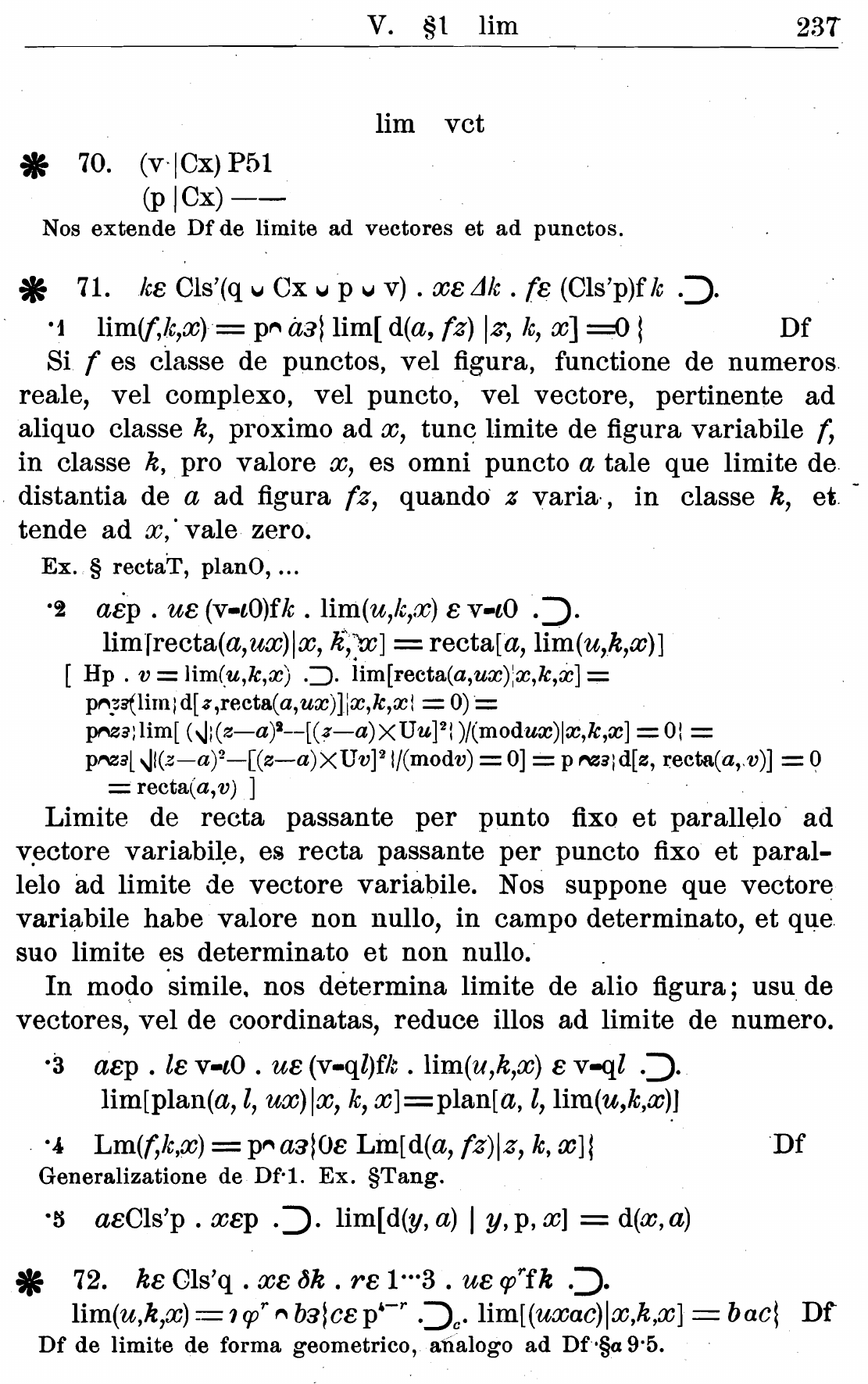}
\caption{FM (p.\,237): lower (n. 71.1) and upper (71.4) limit of sets}
\label{FM237}
\end{center}
\end{figure}

\begin{figure}[htbp]
\begin{center}
\includegraphics[scale=.83]{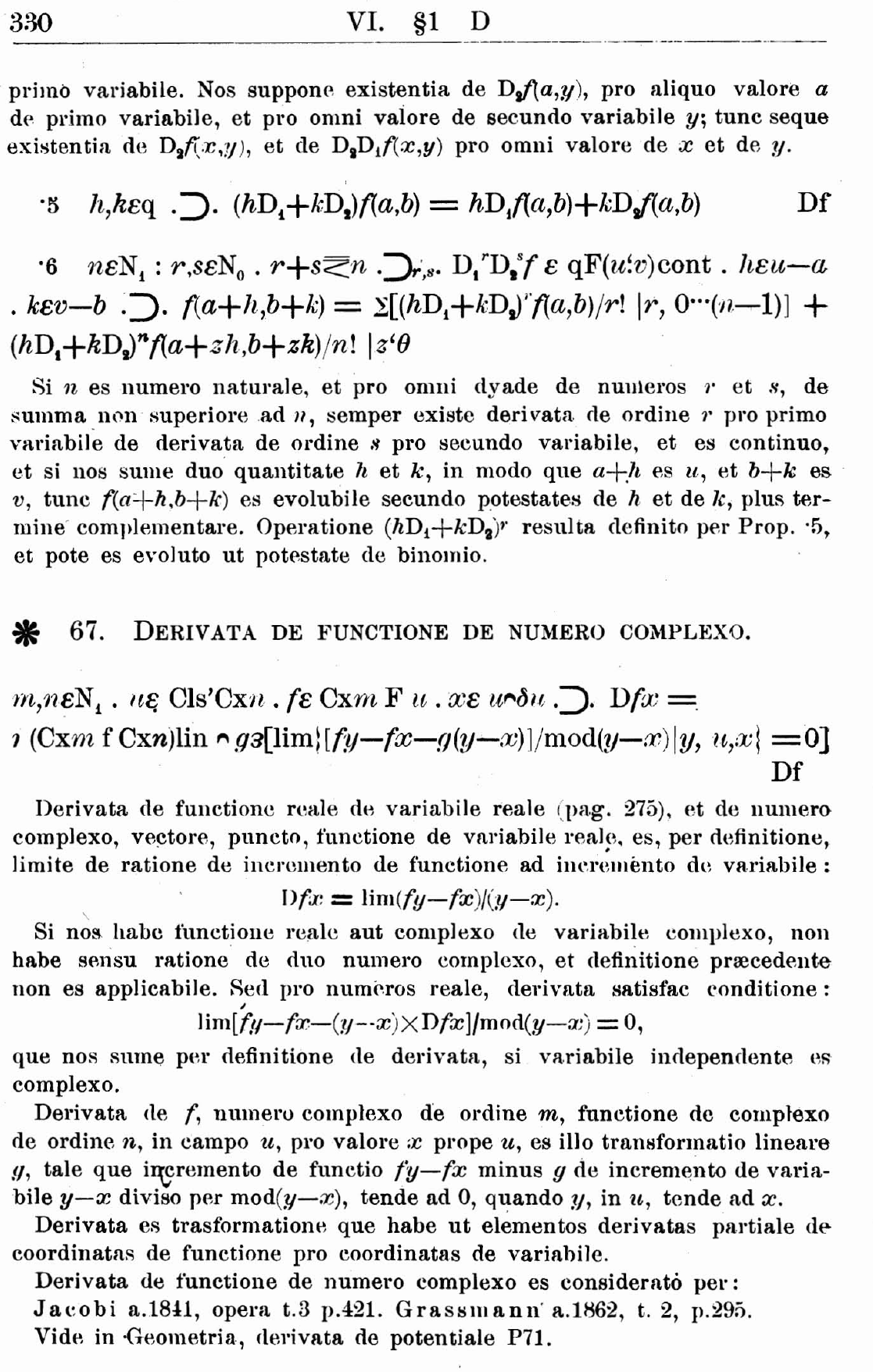}
\caption{FM (p.\,330): Derivative (n.\,67)}
\label{FM330}
\end{center}
\end{figure}

\begin{figure}[htbp]
\begin{center}
\includegraphics[scale=.60]{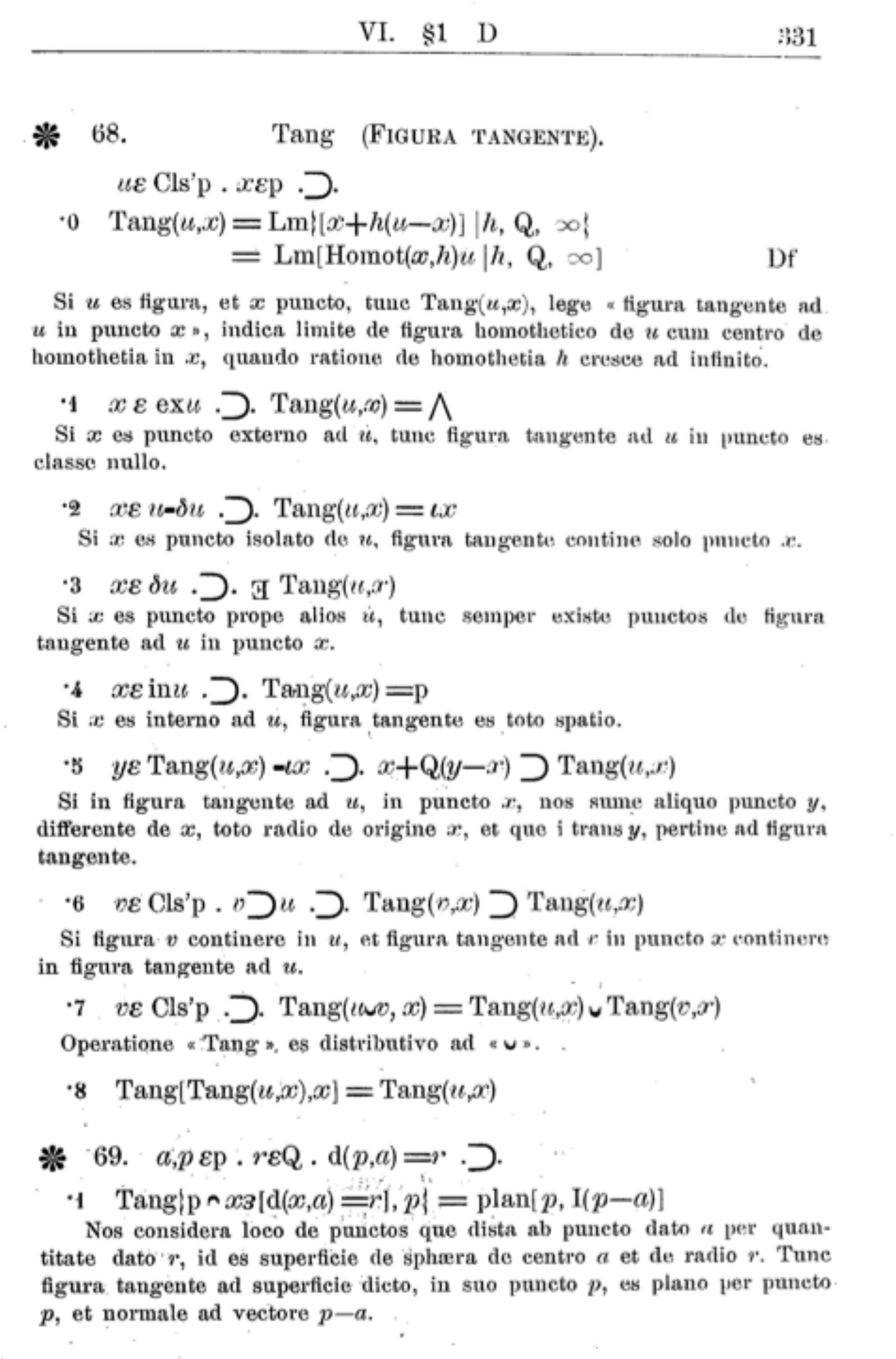}
\caption{FM (p.\,331): Upper tangent affine cone (n.\,68)}
\label{FM331}
\end{center}
\end{figure}

\begin{figure}[htbp]
\begin{center}
\includegraphics[scale=.83]{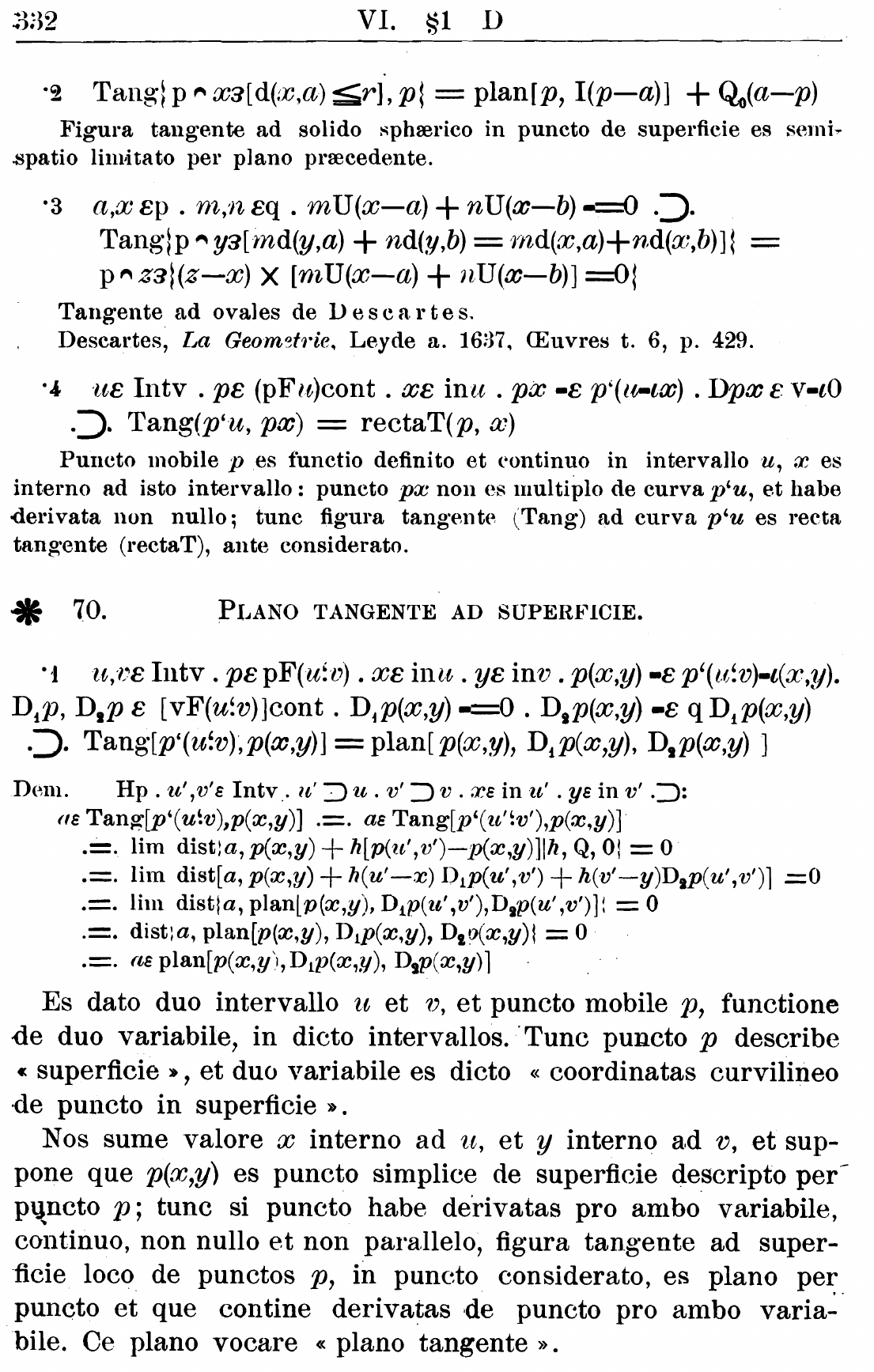}
\caption{FM (p.\,332): Calculus of tangent lines and planes}
\label{FM332}
\end{center}
\end{figure}

\begin{figure}[htbp]
\begin{center}
\includegraphics[scale=.83]{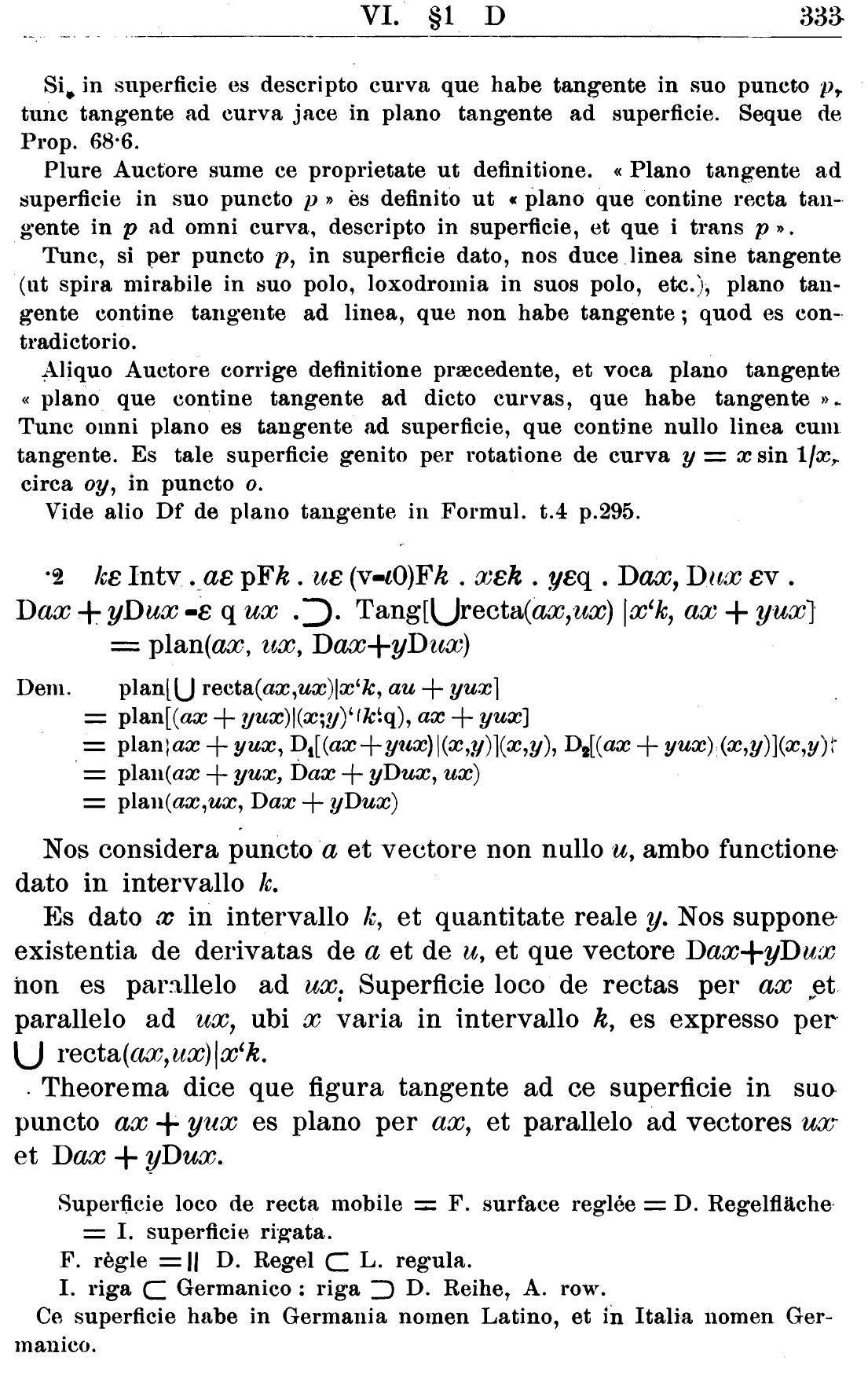}
\caption{FM: Calculus of a tangent plane}
\label{FM333}
\end{center}
\end{figure}

\begin{figure}[htbp]
\begin{center}
\includegraphics[scale=.83]{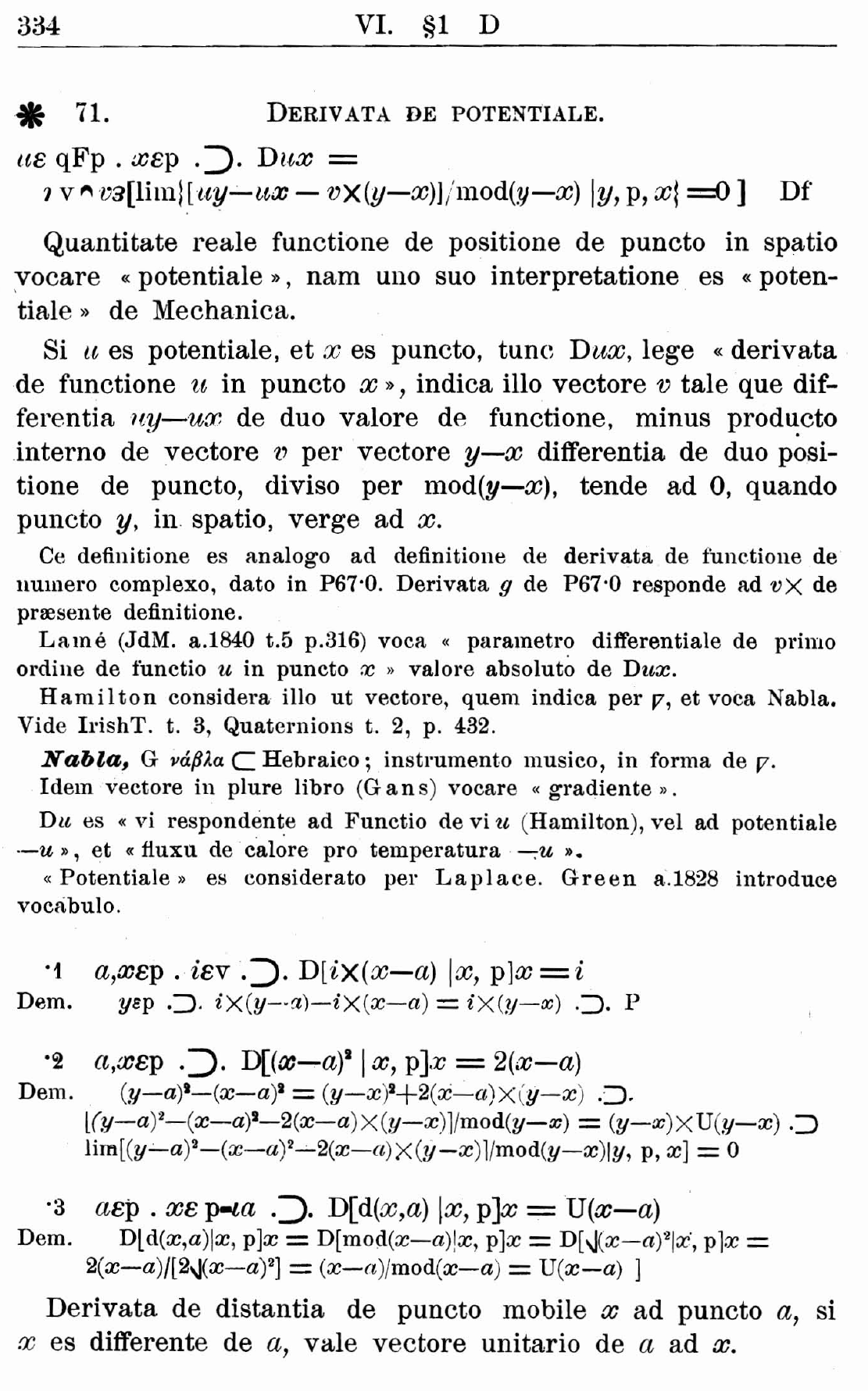}
\caption{FM (p.\,334): Derivative and potential}
\label{FM334}
\end{center}
\end{figure}

\begin{figure}[htbp]
\begin{center}
\includegraphics[scale=.83]{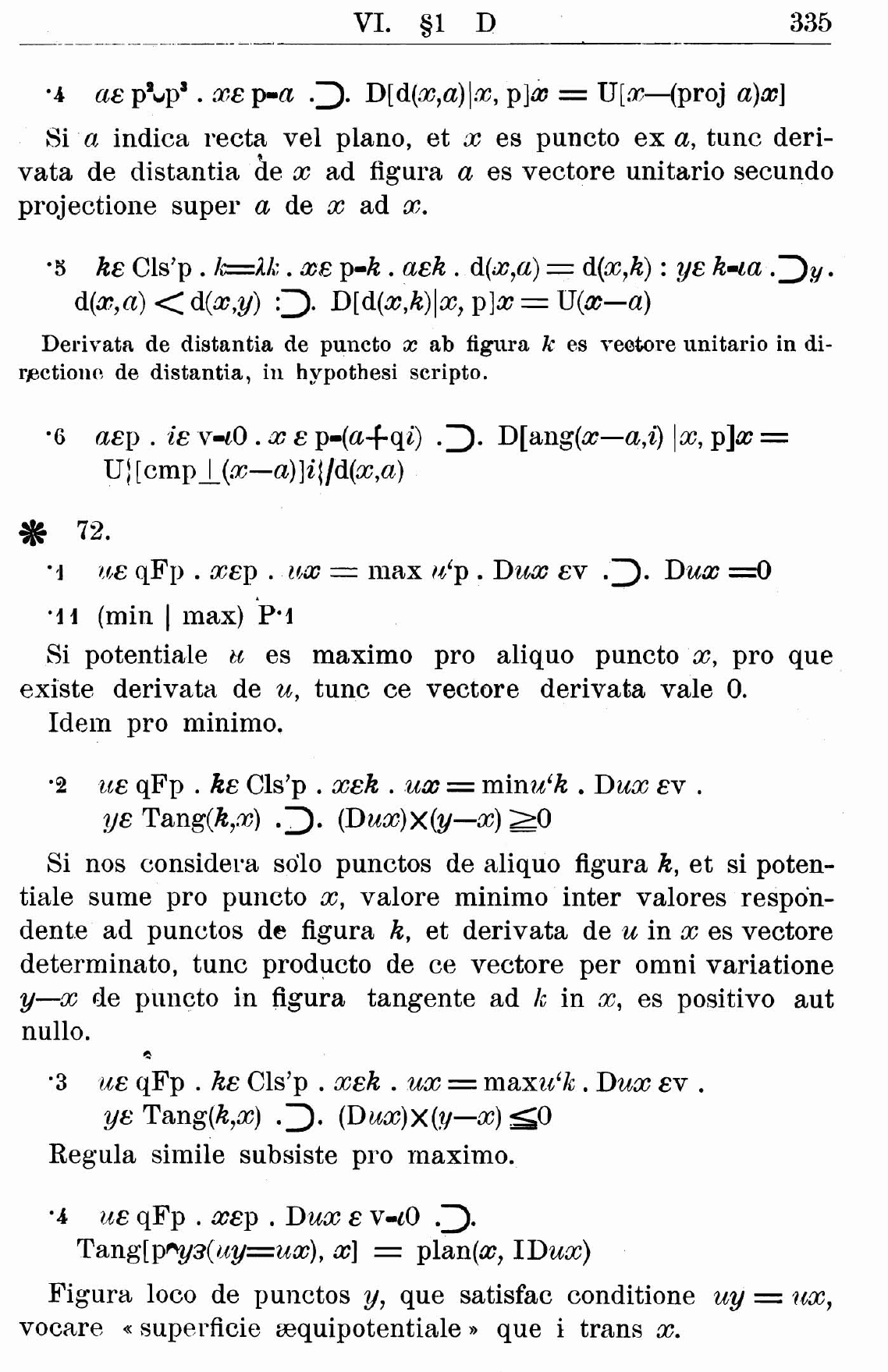}
\caption{FM (p.\,335): Regula n. 72.2 (min) and n. 72.3(max)}
\label{FM335}
\end{center}
\end{figure}

\begin{figure}[htbp]
\begin{center}
\includegraphics[scale=.83]{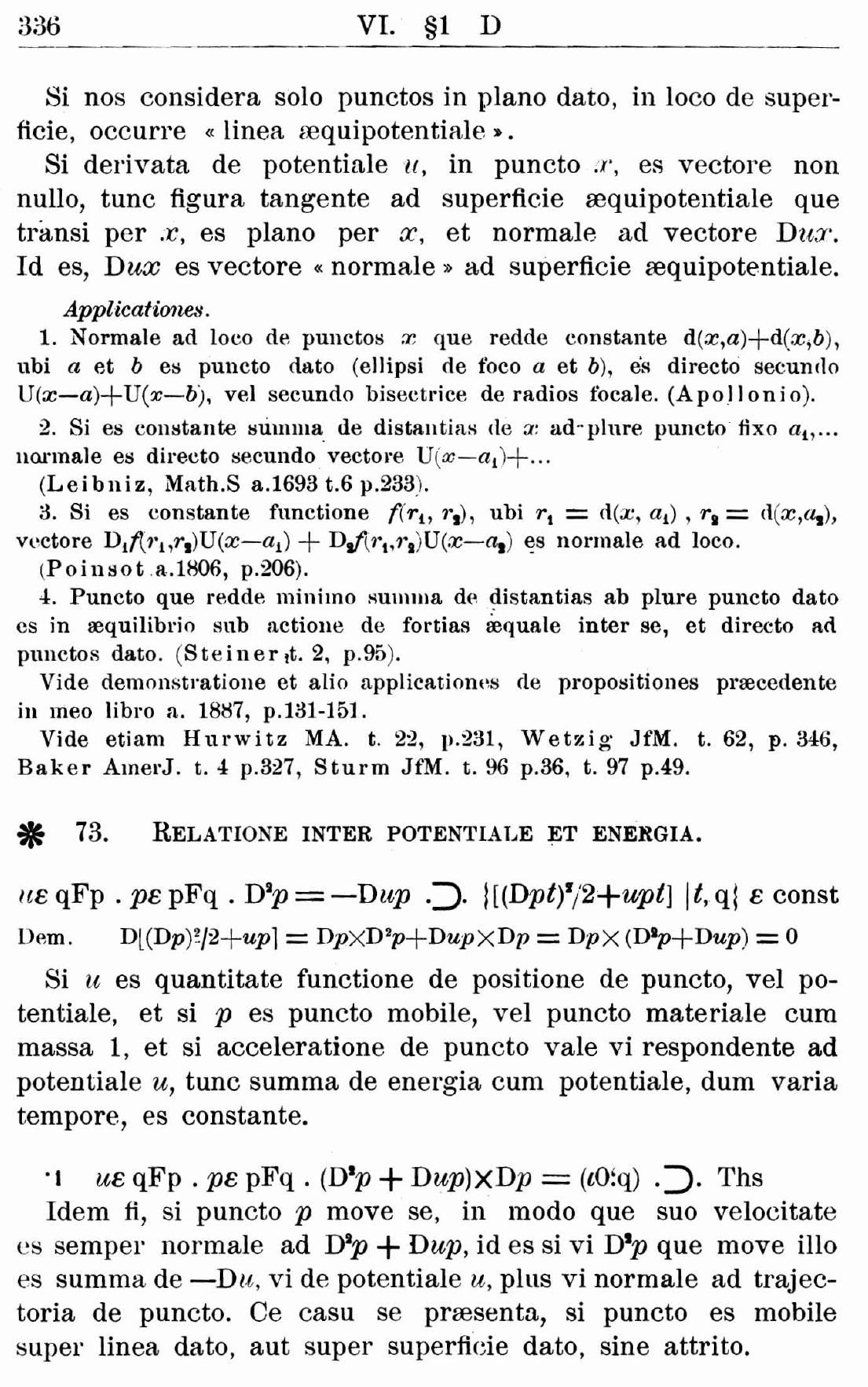}
\caption{FM (p.\,336): Applications}
\label{FM336}
\end{center}
\end{figure}

\end{document}